\documentclass[11pt]{article}

\usepackage{amssymb,color,amsmath,latexsym,eucal,cite}
\usepackage{graphicx,psfrag} 
\usepackage{url}
\usepackage{epstopdf}

\usepackage[ruled]{algorithm}

\newtheorem{ipotesi}{Assumption}[section]    

\newtheorem{theorem}{Theorem}[section]
\newtheorem{lemma}[theorem]{Lemma}

\newtheorem{corollary}[theorem]{Corollary}

\newcommand{\bpr}{{\bf Proof.} \hspace{1.5mm}}
\newcommand{\epr}{\hfill $\Box$ \vspace*{1em}}

\def\IR{\hbox{\rm I\kern-.2em\hbox{\rm R}}}
\def\IN{\hbox{\rm I\kern-.2em\hbox{\rm N}}}

 \def\argmin{\mathop{\rm argmin }}
 \def\liminf{\mathop{\rm liminf }}

\def\t{\theta}
\def\Nk{N_k}
\def\Ntilde{\widetilde N_{k+1}}
\def\Nkp1{N_{k+1}}
\def\Dkp1{D_{k+1}}
\def\pred{{\rm{Pred}}_k}
\def\ared{{\rm{Ared}}_k}
\def\kT{\kappa_T}
\newcommand{\comment}[1]{}

\newcommand{\z}{\phantom{0}}

\newcommand{\be}{\begin{equation}}
\newcommand{\ee}{\end{equation}}

\newcommand{\eqdef}{\stackrel{\rm def}{=}}

\title{Inexact restoration with subsampled trust-region methods for finite-sum minimization\footnotemark[3]
}
\author{ Stefania Bellavia\footnotemark[1], Nata$\check{{\rm s}}$a Kreji\'c\footnotemark[2], Benedetta Morini\footnotemark[1]}
\footnotetext[1]{Dipartimento  di Ingegneria Industriale, Universit\`a degli Studi di Firenze,
Viale G.B. Morgagni 40,  50134 Firenze,  Italia. Members of the INdAM Research Group GNCS. Emails:
stefania.bellavia@unifi.it, benedetta.morini@unifi.it}
\footnotetext[2]{ Department of Mathematics and Informatics, Faculty of Sciences, University of Novi Sad, Trg Dositeja
Obradovi\'ca 4, 21000 Novi Sad, Serbia, Email: natasak@uns.ac.rs. }
\footnotetext[3]{The work of Bellavia and Morini was supported by  Gruppo Nazionale per
il Calcolo Scientifico   (GNCS-INdAM) of Italy. The work of the second author was  supported by Serbian
Ministry of Education, Science and Technological Development, grant no. 174030.

Part of the research was
conducted during a visit of the second author  at Dipartimento di Ingegneria Industriale
supported by Piano di Internazionalizzazione,  Universit\`a degli Studi di Firenze.
}

\begin{document}

%
%
%
%


\maketitle

\begin{abstract}
Convex and nonconvex finite-sum minimization arises in many scientific computing and machine learning applications. 
Recently, first-order and second-order methods where objective functions, gradients  and Hessians are approximated by randomly sampling components of the sum
have received great attention.

We propose a new trust-region method which employs suitable approximations of the objective function, gradient and Hessian   built via random subsampling techniques. The  choice of  the sample size is deterministic and ruled by the inexact restoration approach. 
We discuss  local and global properties for finding approximate first- and second-order optimal points and  function evaluation complexity results.  Numerical experience shows that the new procedure is more efficient,   in terms of  overall computational cost, than the standard trust-region scheme with subsampled Hessians. 
\end{abstract}



\section{Introduction}
The problem we consider in this paper is the following
\begin{equation}\label{minf}
\min_{x\in \IR^n}f_N(x)= \frac{1}{N}\sum_{i=1}^N \phi_i(x),
\end{equation}

\newpage \noindent
where $N$ is very large and finite and $\phi_i: \IR^n\rightarrow \IR$. A number of important problems can be stated in this form, to start with problems in machine learning like classification problems, data fitting problems, sample average approximation of the objective function given in the form of mathematical expectation and so on.

The practical relevance of (\ref{minf}) resulted in a number of methods that are adjusted to this particular form of the objective function. In fact, for very large $ N $ the cost of evaluating $ f _N$ might be really high and the same is true for the gradient and even more for the Hessian evaluation. Therefore a number of methods that use   approximate objective  functions and/or first and second order 
derivatives,  formed  by   partial sums,  is proposed and analysed in literature, see e.g., \cite{bgm,bgmt,bkkj,bbn_2017, 
bkm2018,bcms,bbn,nocedalsurvey,noc3,noc2,em,Friedlander,llh,rm,xu1,XuRoosMaho17}.

Concerning the approximation of  the objective function,  one of the possible approaches  is to use relatively rough approximations  at early stages of the optimization procedure and gradually increase  the accuracy   to arrive at full precision at the late stage  of the iterative procedure;
the gradient is approximated accordingly.
This way one hopes to save computational effort and yet to solve the original problem eventually.  
Very often the term scheduling is used to describe  the approximation of the objective function by means of  a partial sum. There is a number of algorithms proposed for the scheduling problem, ranging from simple heuristics that increase the number of terms in the partial sum that approximates the objective function by a certain percentage in each iteration, \cite{bkkj,bbn,Friedlander,rm} to more elaborate schemes that connect the progress achieved during the optimization procedure to the number of terms in the partial sum \cite{B.Tuan, B.teorija, bgm, bgmt, bkkj, bkm2016, bkm2018, bcms, noc2, Feris, nas, nas1, nkjmm, P2,Polak}. 


Besides the problem of scheduling, one has to decide between first- and second-order optimization method to be employed. A detailed survey is presented in \cite{nocedalsurvey}. 
A number of first-order methods has been proposed and analysed in the literature.
Given that the main cost comes from large $ N $ one might be tempted to conclude that computing Hessians, or some other second order information might be prohibitively costly  and thus opt for a first order method, especially if the problem (\ref{minf}) should be solved with limited precision.  However, recently there has been reported in several papers that careful adjustment and implementation of second order methods might be worth considering if the true Hessian is approximated by a partial sum of Hessians  $ \nabla^2 \phi_i(x) $ consisting of a  significantly smaller number of terms than $ N$.  
This way one can generate useful information with significantly smaller cost than the true Hessian and get enough advantage over   first-order methods in terms of resilience to problem ill-conditioning and low sensitivity to parameter tuning, \cite{bbn_2017,bkkj, bbn,noc2, noc3,em,Pilanci,xu1,XuRoosMaho17,rm}. 

The method we present here combines the Inexact Restoration (IR) framework with the trust-region optimization method  \cite{cgt} to 
simultaneously design the scheduling and the optimization procedure for solving (\ref{minf}) and  represents a new
approach for the problem under consideration.

The Inexact Restoration  method,   introduced in \cite{jmmpillota}, 
is a constrained optimization tool particularly suitable for problems where one does not want to enforce feasibility in all iterations. 
The key idea of the IR approach is to treat feasibility and optimality in a modular way and to improve each one in  separate procedures;
the combination of feasibility and optimality is then monitored through a suitable  merit function. 
Each iteration ensures the sufficient decrease of a suitable merit function and therefore, under certain assumption, convergence to a feasible optimal point.
In \cite{jmm,jmmpillota} the combination of the IR strategy  with  trust-region methods is proposed and analysed   for general constrained problems.

The  application of IR strategy to the unconstrained optimization  problem (\ref{minf}) requires its reformulation as a constrained problem.
Letting $I_M$ be an arbitrary nonempty  subset of   $\{1, \ldots, N\}$ of cardinality $| I_M|$ equal to $M$, we 
reformulate problem (\ref{minf}) as 
\begin{equation} 
\begin{aligned}
 & \min_{x\in \IR^n}  f_M(x)   =    \frac{1}{M}\sum_{i\in I_M}  \phi_i(x).\\
 & \mbox{ s.t. }      M=N,    
\end{aligned}  \label{minf1}
\end{equation}
 Evaluating infeasibility in (\ref{minf1}) is cheap while computing the objective function is   expensive whenever $M$ is large.
Thus, using the reasoning from \cite{jmm,jmmpillota} we define a new algorithm  that exploits the structure of the problem considered and  
takes advantage of the  modular structure of IR and the trust-region  optimization method at the same time.  
 Specifically, the trust-region mechanism is applied  to model  $f_M$ at each iteration and the 
IR framework is applied to test for  the acceptance  of the iterates and to determine  
the scheduling sequence, i.e. the value of $M$ through the iterations.  
The test acceptance of the new iterate allows us to deal with inaccuracy in function and derivatives. In particular,
the number of terms in the partial sum is fixed at the beginning of each iteration  in the restoration   phase and possibly changed in the optimality phase  where the trial iterate is computed.

Clearly, the higher feasibility is the more accurate   $f_M$ is with respect to $f_N$.  The new procedure has two important properties:  
partial sums,  possibly consisting of small sets of $ \phi_i$'s, can be used in the early stage of the iterative procedure to decrease the computational cost;  
the original objective function in (\ref{minf}) is recovered for all iteration indices large enough, thus  allowing for 
the solution of the given problem. Clearly, when full precision of the objective function and the gradient is reached, one can rely on 
the theory and machinery of standard trust-region methods \cite{cgt}. 

The scheme presented here applies to both first- and second-order trust-region models.
If a linear model is used, the resulting procedure is a subsampled gradient method with variable stepsize.  
When second-order models are used, the Hessian can be approximated using a subset of the sample used to approximate function and gradient.
The error in such Hessian approximation plays an important role in the asymptotic convergence rate. 
In the case of strongly  convex problems,    the analysis  for local  linear  convergence rate
is presented, both in deterministic and probabilistic settings, and  an adaptive choice of the sample for  
Hessian approximation is proposed.

We also  provide a  function evaluation  complexity result which  resembles the classical result for the trust-region methods for (\ref{minf}) and the results obtained in \cite{bkm2018}.
It is shown that at most  $O( \varepsilon^{-2}) $ evaluations of the  possibly subsampled function $f_M$, $M\le N$, and its derivatives are needed to compute a first-order approximate critical point. Then the worst-case complexity of the standard trust-region is recovered with expected significant computational savings due to scheduling.

 Our approach considerably differs from the IR procedure and trust-region method in  \cite{jmm,jmmpillota} since the objective function in our formulation changes with $M$ through the iterations.  It also differs from IR approaches in  \cite{nkjmm, bkm2016, bkm2018} that employ approximate  objective function and its derivatives
and have been successfully applied to  constrained and unconstrained  problems,  including  problem (\ref{minf}); in papers \cite{nkjmm, bkm2016} the IR is combined with  a line search strategy, while in  \cite{bkm2018} the considered problem is constrained and regularization techniques are used in the optimization phase. The approach presented here relays on \cite{bkm2018} in terms of general idea but the problem is more specific being a finite-sum  rather than a general objective function computed approximately and being unconstrained. These specifications allow us to design an efficient sample update rule which is connected with the trust-region size.

  The value of $M$ is fixed via a deterministic rule while  the trust-region schemes 
in \cite{llh, XuRoosMaho17, bcms},  approximating either  functions, gradients and Hessians   \cite{llh,bcms} or  Hessians only  \cite{XuRoosMaho17}, are designed using sample sets whose cardinality is determined by high probability and nonasymptotic convergence analysis. 

The nature of IR allows changes in  the feasibility  through iterations and the change is not necessarily monotone, 
i.e., the cardinality of the  subset that defines the approximate objective can both increase and decrease, 
depending on the feedback from the trust-region progress  made in each iteration.
The case where $M$ is increased by a prefixed percentage at each iteration is a particular case of our strategy.   
 In this latter case  our method differs from a straightforward subsampled trust-region procedure with increasing sample size   
in both the merit function and the acceptance criterion. Remarkably, their employment  allow to prove  optimal complexity results that otherwise require adaptive accuracy requirements \cite{bcms}.

 
This paper is organized as follows. In Section 2 we present our method and prove that it  is well defined.
Furthermore, we prove that full accuracy is eventually reached and that the set of standard assumptions yield first-order stationary points. Some  issues concerning the realization of the procedure are considered in Section 3; the scheduling rule is modified to avoid unproductive decrease  in precision and a discussion on first and second order trust-region models is provided. Section 4 deals with strongly convex problems; 
we prove q-linear convergence as well as q-linear convergence in expectation 
under   probabilistic bounds for Hessian subsampling. Section 5 provides worst-case function evaluation complexity.
The numerical performance of the proposed method is tested on a set of classification problems and the results are reported in Section 6.

\section{The Algorithm}
Let $I_M$ be an arbitrary nonempty  subset of   $\{1, \ldots, N\}$ of cardinality $| I_M|$ equal to $M,$
$$
 I_M\subseteq\{1, \ldots, N\}, \quad|I_M|=M, \ \ M\ge 1,
$$
and reformulate (\ref{minf}) as the constrained problem (\ref{minf1}).
We measure the level of infeasibility with respect to the constraint $M=N$ by the function $h$  with the following properties.
\begin{ipotesi}\label{assh}
Let $h:\IN\rightarrow \IR$ be   a monotone, strictly decreasing function such that $h(1)>0$, $h(N)=0$.   
\end{ipotesi}
This  assumption implies  
\begin{equation}\label{boundh}
\underline{h}\le h(M) \ \ \mbox{ if } \ \ 0<M<N, \quad \mbox { and }\quad  h(M)   \le \bar h \ \ \mbox{ if } \ \ 0<M\le N,
\end{equation}
for $M\in \IN$ and  $\underline{h}=h(N-1)$ and $\bar h=h(1)$.
One possible choice for $h$ is $h(M)=(N-M)/N, \   0<M\le N$. 

 Suppose $\phi_i$, $1\le i\le N$, be continuously differentiable and let  
$\|\cdot\|$ denote  the 2-norm.

The method introduced in this section combines the Inexact Restoration,  an  approach for optimization of functions evaluated inexactly,
with the trust-region methods. We will refer to it as {\sc iretr}. It employs the merit function 
\begin{equation}\label{merit}
 \Psi(x,M,\t)=\t f_M(x)+(1-\t)h(M),
\end{equation}
with $\t\in (0,1)$ and aims to minimize both $f_M$ and the infeasibility $h$. 
Since the reductions in the values of $f_M$ and $h$ may not be achieved  simultaneously, a weight $\t$ is used and 
a trust-region method is employed  to generate
a sequence $\{(x_k, \Nk, \t_k)\}$ such  that $ \Psi(x_k, \Nk, \t_k)<\Psi(x_{k-1}, N_{k-1}, \t_k)$.
The main theoretical properties of the new method, shown in the next section,  are:  the sequence $\{\t_k\}$ is nonicreasing and uniformly bounded away from zero, 
$N_k=N$ for all $k$ sufficiently large and $\|\nabla f_N(x_k)\|\rightarrow 0$ as $k\rightarrow \infty$.

Concerning the trust-region problem,  suppose that   $x_k$ is given. Then,    a trial sample size $\Nkp1$ is selected,   
$I_{\Nkp1}\subseteq\{1, \ldots, N\}$ is chosen and the model 
$m_k(p)$  for $f_{\Nkp1}$ around $x_k$  of the form
\begin{equation}\label{model}
m_k(p)=f_{\Nkp1}(x_k)+\nabla f_{\Nkp1}(x_k)^Tp+\frac 1 2 p^T B_{k+1} p,
\end{equation}
is built. Here $\nabla f_{\Nkp1}$ denotes the gradient of $f_{\Nkp1}$ and  $B_{k+1}\in \IR^{n\times n}$ is a symmetric 
approximation to the Hessian  $\nabla^2 f_{\Nkp1}(x_k)$  in case $\phi_i$, $1\le i\le N$, are  twice  continuously differentiable.
Trivially $m_k(0)=f_{\Nkp1}(x_k)$ and the smaller $h(\Nkp1)$, the larger  becomes the accuracy in the approximation to $f_N$ and $\nabla f_N$.
Then, letting $\Delta_k>0$ denote the trust-region radius and ${\cal B}_k=\{x_k+p\in \IR^n: \|p\|\le \Delta_k\}$ be  the trust-region,
the trust-region  problem is 
\begin{eqnarray}\label{tr_pb}
\min_{ \|p\|\le \Delta_k} m_k(p).
\end{eqnarray}
As in the standard trust-region schemes, problem (\ref{tr_pb}) is solved approximately and 
the computed step $p_k$ is required to provide a sufficient reduction 
in the model in terms of  the Cauchy step $p_k^C$, i.e., the minimizer of  the model $m_k$ 
along the steepest descent  $-\nabla f_{\Nkp1}(x_k) $ within ${\cal B}_k$
\begin{eqnarray}\label{cauchy}
p_k^C=\argmin_{\begin{small}
\begin{array}{c} p= -t \nabla f_{{\Nkp1}}(x_k), \, t>0\\ \|p\|\le    \Delta_k \end{array}
\end{small}} m_k(p).
\end{eqnarray}
 Then, if a sufficient reduction in the function $\Psi$ is achieved, the step $p_k$ is accepted and the new iterate
$x_{k+1}$ is set equal to $x_k+p_k$.  Otherwise,  the step is rejected and the trust-region radius is reduced. The  specific form of the predicted and actual reduction used in the acceptance criterion will be given below, after detailing the Algorithm's steps. 

Now we present the new Algorithm {\sc iretr} which  aims at finding an
$\varepsilon_g$--accurate first-order optimality point defined as follows
\be \label{epsfo}
 \|\nabla f_{\Nkp1}(x_k)\|\le \varepsilon_g  \quad \mbox { and } \quad \Nk=N,
\ee 
and comment on it, see Algorithm \ref{IRTR_algo}. 

 \begin{algorithm}
\caption{The  algorithm {\sc iretr}}\label{IRTR_algo}


 \vskip 5 pt 
 \noindent
 Given   $x_0\in \IR^n$, $N_0$ integer in $(0,N]$, $\t_0\in (0,1)$,  $B_0\in \IR^{n\times n}$, 
 $\Delta_0>0$, \\
$0<\zeta_1<1<\zeta_2$, $\gamma\in (0,1],\,  r, \eta, \tau \in (0,1)$,  $\mu \in [0,1)$ the accuracy level $\varepsilon_g\ge 0$.
\vskip 4pt
\noindent
0. Set $k=0$, ${\cal T}_0=0$, $\Delta_0=\Delta_0^{({\cal T}_0)}$;\\
1.  If $\Nk<N$, find $\Ntilde$ such that $\Nk<\Ntilde\le N$, and 
\begin{equation}\label{feas}
h(\Ntilde)\le r h(\Nk).
\end{equation}
\hspace*{11pt}   If $\Nk=N$, set $\Ntilde =N$.  
\vskip 2pt \noindent
 2.    Choose $\Nkp1$ such that $\Nkp1\le \Ntilde$, and 
\begin{eqnarray}
 h(\Nkp1)-h(\Ntilde) &\le& \mu  \left(\Delta_k^{({\cal T}_k)}\right)^{1+\gamma}.\label{new2}   
\end{eqnarray}
\hspace*{11pt}   If $\Nk=N$ and $\|\nabla f_{\Nkp1}(x_k)\|\le \varepsilon_g$, stop.
\\
\hspace*{11pt} Build the  model $m_k(p)$ for $f_{\Nkp1}(x_k) $ in (\ref{model}).\\
\hspace*{11pt} Find an approximate trust-region solution $p_k$ such that
\begin{eqnarray}
m_k(0)-m_k(p_k) &\ge& \tau (m_k(0)-m_k(p_k^C)) \label{steptr} 
\end{eqnarray}
\hspace*{13pt}where $p_k^C$ is given in \eqref{cauchy}.\\
3. If $\Nk=N$ and $\Nkp1<N$  and
\begin{eqnarray}
f_{N}(x_k)-m_k(p_k) &<& \tau (m_k(0)-m_k(p_k^C))  \label{new1} 
\end{eqnarray}
\hspace*{15pt} take  $\Delta_k^{({\cal T}_k+1)}=\zeta_1\Delta_k^{({\cal T}_k)}$, set ${\cal T}_k={\cal T}_k+1$ and repeat Step 2.    
\\
4.  Compute the penalty parameter $ \theta_{k+1}$
\begin{equation}
\begin{aligned}
\t_{k+1}=
\left\{
\begin{array}{ll}
&\t_k  \hspace*{30pt} \mbox{  if } \  \pred(\t_k)\ge\eta( h(\Nk)-h(\Ntilde))\\
& \displaystyle \frac{(1-\eta)(h(\Nk)-h(\Ntilde))}{m_k(p_k)-f_{\Nk}(x_k)+h(\Nk)-h(\Ntilde)} \quad  \mbox{otherwise}.
\end{array}
\right.
\end{aligned}\label{tkp1}
\end{equation}
5.  If 
\begin{equation}\label{ared}
\ared(\t_{k+1})\ge \eta \pred(\t_{k+1}),
\end{equation}
\hspace*{14 pt}Set $x_{k+1}=x_k+p_k$, $\Delta_k=\Delta_k^{({\cal T}_k)}$.  Choose  $ \Delta_{k+1}^{(0)}\in [\Delta_k,\zeta_2\Delta_k ]$, \\
\hspace*{14 pt}set $k=k+1,{\cal T}_k=0$, and go to Step 1.\\
\hspace*{10pt}
Else take  $\Delta_k^{({\cal T}_k+1)}=\zeta_1\Delta_k^{({\cal T}_k)}$, set ${\cal T}_k={\cal T}_k+1$ 
and go to Step 2. 
\end{algorithm}

\vskip 10pt
Given  $x_k$, $N_k$ and $\t_k$  we describe the  $k$th iteration. In Step 1 the feasibility is improved. If $N_k<N$, 
we predict the cardinality $\Ntilde$ such that  the value  $h(\Ntilde)$ is smaller than $h(\Nk)$ and  
at most equal to a prefixed fraction of $h(\Nk)$. 
 In case  $h(M)=(N-M)/N, \   0<M\le N$, taking into account that  $N_k$ and $\tilde N_{k+1}$ are integers
it can be shown that condition (\ref{feas}) holds if and only if $0< N_k<\tilde N_{k+1}$ provided that 
$ h(2)/h(1) < r <1$.

In Step 2, an attempt is made to reduce the computational effort i.e. to enlarge infesibility;  $\Nkp1$ is chosen such that $\Nkp1\le \Ntilde$ and the  bounded
deterioration (\ref{new2}) on the value of $h(\Nkp1)$ with respect to $h(\Ntilde)$ is imposed. 
In principal such  control allows  us to reduce $\Nkp1$ below both $\Nk$ and $\Ntilde$. 
On the other hand, the upper bound in (\ref{new2}) 
depends on the trust-region radius and $\Nkp1$ will be equal  to $\Ntilde$ whenever $\Delta_k$ is small enough.
 If $N_k=N$, the stopping criterion $\|\nabla f_{\Nkp1}(x_k)\|\le \varepsilon_g$ is checked. 
This is  supported by the fact that, when  $N_k=N$, we may expect $\Nkp1$ 
be close to $N$ and $ \nabla f_{\Nkp1}(x_k)$ be close to  $ \nabla f_{N}(x_k)$ in a probabilistic sense; we will further discuss this issue in Section  \ref{Sec3}.
If  (\ref{epsfo}) is not met, using $I_{\Nkp1}\subseteq\{1, \ldots, N\}$, the trust-region model $m_k(p)$ is built and  (\ref{tr_pb}) is approximately solved.
The computed step $p_k$ is required to provide the sufficient reduction (\ref{steptr}) 
in the model in terms of  the Cauchy step $p_k^C$. 

   The acceptance rule for $p_k$  in Step 5 depends on the predicted and actual reduction defined as follows:
\begin{eqnarray}
\pred(\t)&=&\t(f_{\Nk}(x_k)-m_k(p_k))+(1-\t)(h(\Nk)-h(\Ntilde)), \label{defpred}\\
\ared(\t) &=&  \Psi(x_k, \Nk, \t )- \Psi(x_{k}+p_k, \Nkp1, \t ) \nonumber\\
&=&  \t(f_{\Nk}(x_k)-f_{\Nkp1}(x_k+p_k))+(1-\t)(h(\Nk)-h(\Nkp1)),\label{defared}
\end{eqnarray}
where the last equality follows from \eqref{merit}. 
We observe that $\pred$ uses the last accepted values $f_{\Nk}(x_k)$ and $\Nk$
and is a linear combination of two predicted values:  the predicted model decrease $ f_{\Nk}(x_k)-m_k(p_k)$
and the predicted   infeasibility  decrease $h(\Nk)-h(\Ntilde)$. As for $\ared$, given $\theta$, it measures the actual reduction of $\Psi$.

The new penalty parameter $ \t_{k+1} $ computed in Step 4 is the largest value  that ensures 
\be \label{pred} 
\pred(\t_{k+1}) \geq \eta(h(N_k) - h(\widetilde{N}_{k+1}))\ge 0,
\ee
as $h(N_k) - h(\widetilde{N}_{k+1})\ge 0$ by (\ref{feas}).  In case $N_k<\Ntilde$ such condition implies   $\pred(\t_{k+1})$  strictly positive.  In case
$N_k= \Ntilde=N$, $\pred(\t)$ reduces to $\t(f_{N}(x_k)-m_k(p_k))$ and from \eqref{steptr}  it follows  $\pred(\t)\ge \tau \t(m_k(0)-m_k(p_k^C))>0$  whenever $\Nkp1=N$. On the other hand, in case  $N_k=\Ntilde=N$ and $\Nkp1<N$,
Step 3 is necessary to enforce positivity of $\pred(\t_{k+1})$ as  $m_k(0)=f_{\Nkp1}(x_k)\ne f_N(x_k)$.
 In fact, $\pred(\t)>0$    follows from taking a step such that 
$f_{N}(x_k)-m_k(p_k) \ge  \tau (m_k(0)-m_k(p_k^C)) $.
We further notice that attempting 
$\Nkp1<N$ when $N_k=N$ is meaningful if  the model is a good approximation of $ f_N $ around $x_k$  and thus one can expect some progress, or at least a limited  deterioration in the value of the full  objective function $f_N$.  Enforcing $f_{N}(x_k)-m_k(p_k) \ge  \tau (m_k(0)-m_k(p_k^C)) $ is a minimal requirement  on the agreement between  $f_N$ at $x_k$ and the model at the trial step.

Finally, in Step 5 the step $p_k$ is accepted  if the ratio between the predicted  reduction $\pred(\t_{k+1}) $ and the actual reduction $Ared_k(\t_{k+1}) $ is larger than a prefixed scalar $\eta$,
otherwise the trust-region radius is reduced and the procedure is repeated starting from Step 2. 

 Notice that the trust-region size can be reduced several times during one iteration, i.e., only successful iterations yield to the increment of the iteration counter $k$. To emphasize this fact,  within each iteration, we introduce an additional counter $ {\cal T}_k $ for the number of decreases of the trust-region size. The feasibility measure $ N_{k+1} $  might be modified  several times within one iteration as well, but changes due to (\ref{new2}) and (\ref{new1}) do not necessarily correspond  to the number of reductions of the trust-region size.   The penalty parameter $ \theta_k$ has an analogous behaviour. 
For this reason and to avoid notation clustering, we  do not introduce additional counters for $ N_{k+1} $ and $ \theta_{k+1} $ within the same iteration. 

We start the analysis of the new method proving that the $k$th iteration of  Algorithm {\sc iretr}  is well defined  since appropriate values of $ N_{k+1} $ and $ \theta_{k+1} $ will be reached in a finite number of attempts.
Here and in Section 5,  $B_{k+1}$ can  be the null matrix and our  analysis covers the use of both  first-order and second-order models.
 
\begin{lemma}\label{ChoiceNkp1}
Steps 2 and 3 of Algorithm {\sc iretr} are well-defined.
\end{lemma}
\bpr
For any positive $\Delta_k^{({\cal T}_k)}$ inequality  (\ref{new2}) trivially holds in the limit case
$\Nkp1=\Ntilde$.  Analogously, Step 3 can not be repeated infinitely many times as  for $ {\cal T}_k $ large enough, 
$ \Delta_k^{({\cal T}_k )}  $ will be small enough to yield $\Nkp1=\Ntilde=N$.
\epr

We now make the following assumption.
\begin{ipotesi}\label{assxk}
$\{x_k\}\subset \Omega$ where $\Omega$ is a compact set in $\IR^n$.
\end{ipotesi}

\begin{lemma}\label{lemmatheta} 
Let Assumptions \ref{assh} and \ref{assxk} hold. Suppose that $ \phi_i$, $1\le i\le N$, are  continuous  in $\Omega. $ Then the sequence $ \{\theta_k\} $ built  in Algorithm {\sc iretr}  is positive, nonincreasing and bounded  away from zero, $\t_{k+1} \geq\underline{\t}>0$ with $\underline{\t}$ independent of $k$ and 
\eqref{pred} holds. 

\end{lemma}
\bpr  We have $\t_0>0$ and proceed  by induction  assuming  that $\t_k$ is positive. First consider the case where $\Nk=\Ntilde$ (equivalently $\Nk=\Ntilde=N$).  Then $h(\Nk)-h(\Ntilde)=0$ and, due to   Step 3,
$\pred(\t)=\t(f_{\Nk}(x_k)-m_k(p_k))>0$   for any positive  $\t$. Thus  $\t_{k+1}=\t_k $ and (\ref{pred}) holds. 

Let now  suppose $\Nk<\Ntilde$. If inequality $\pred(\t_k) \ge  \eta(h(\Nk)-h(\Ntilde))$ holds  then $\theta_{k+1}=\theta_k$
satisfies (\ref{pred}).  Otherwise,  we have 
$$
\t_k  ( f_{\Nk}(x_k)-m_k(p_k)-( h(\Nk)-h(\Ntilde) ))   < {(\eta-1)(h(\Nk)-h(\Ntilde))} ,
$$
and since the right hand-side   is negative by construction, it follows
$$
f_{\Nk}(x_k)-m_k(p_k)-(h(\Nk)-h(\Ntilde))<0.
$$ 
Consequently, $\pred(\t )\ge \eta(h(\Nk)-h(\Ntilde))$ is satisfied if 
$$
\t (f_{\Nk}(x_k)-m_k(p_k)-(h(\Nk)-h(\Ntilde)))\ge  (\eta-1)(h(\Nk)-h(\Ntilde)),
$$
i.e., if 
$$ \t \le \t_{k+1}\eqdef \frac{(1-\eta)(h(\Nk)-h(\Ntilde))}{m_k(p_k)-f_{\Nk}(x_k)+h(\Nk)-h(\Ntilde)}. $$
Hence $\t_{k+1}$ is the largest value satisfying (\ref{pred}) and $ \t_{k+1} < \t_k. $

Let us now prove that $ \t_{k+1} \geq \underline{\t}. $   Using Assumptions  \ref{assxk} and  continuity of $ \phi_i$, $1\le i\le N$,  let
\begin{equation}\label{kappaphi}
\kappa_\phi=   \max_{\begin{small} \begin{array}{c} 1\le i\le N\\ x\in \Omega\end{array} \end{small}}  |\phi_i(x)|.
\end{equation}
Then, using (\ref{boundh}),  for $M$ such that $0<M\le N$ there holds 
\begin{eqnarray*}
f_N(x_k)-f_M(x_k) & = & \frac{1}{N} \sum_{i \in I_N} \phi_{i}(x_k) -  \frac{1}{M} \sum_{i \in I_M} \phi_{i}(x_k) \\
& = & \left( \frac{1}{N} - \frac{1}{M}\right) \sum_{i \in I_M} \phi_{i}(x_k)  + \ \frac{1}{N}  \sum_{i \in I_N \diagdown I_M} \phi_{i}(x_k), 
\end{eqnarray*}
and therefore for any integer $M$, $0<M\le N$
\begin{eqnarray} 
|f_N(x_k)-f_M(x_k)|  & \leq  & \frac{N-M}{NM} M \kappa_\phi + \frac{N-M}{N} \kappa_\phi \nonumber \\
&= & \frac{2(N-M)\kappa_\phi}{ N\,h(M)}  h(M)  \nonumber \\
&\le &  \frac{2(N-M)\kappa_\phi}{ N\, \underline{h}} h(M)  \nonumber  \\
&\le &  \frac{2(N-1)\kappa_\phi}{ N\, \underline{h}} h(M)  \nonumber  \\
&\eqdef & \sigma h(M). \label{newassf}
\end{eqnarray}
Also note that by (\ref{feas}) and (\ref{boundh})  

\begin{equation}\label{riduh}
h(\Nk)-h(\Ntilde)\ge (1-r)h(\Nk)\ge (1-r) \underline{h}.
\end{equation} 
  Moreover, 
\begin{eqnarray*}
 &  & m_k(p_k)-f_{\Nk}(x_k)+h(\Nk)-h(\Ntilde)  \le m_k(p_k)-f_{\Nk}(x_k)+h(\Nk)\\
&  & \le m_k(0)-f_{\Nk}(x_k)+\bar h =  f_{\Nkp1}(x_k)-f_{\Nk}(x_k)+ \bar h  \\ 
&  & \le  |f_{\Nkp1}(x_k)-f_N(x_k)|+|f_N(x_k)-f_{\Nk}(x_k)| +  \bar h\\
&  &  \le  \sigma(h(\Nk)+h(\Nkp1))+  \bar h \le  (2\sigma+1) \bar h,
\end{eqnarray*}
and  $\t_{k+1}$ in (\ref{tkp1}) satisfies 
$$
 \t_{k+1} \ge \frac{(1-\eta)(1-r) \underline{h}}{(2\sigma+1)\bar h}  \eqdef \underline{\theta},
$$
and the proof is completed.
\epr

To establish the well-definiteness of Steps 4 and 5, we make the following assumptions.

\begin{ipotesi}\label{assf}
The gradients  $ \nabla \phi_i$, $1\le i\le N$,  are  Lipschitz continuous on the segments $[x_k, x_k+p_k]$,
for all $ k\ge 0$ and for all $ p_k$ generated in the repetition of Steps 2--5.
\end{ipotesi}

\begin{ipotesi}\label{Bk} There exists positive $\kappa_B$ such that for all $ k $ 
$$\| B_{k+1}\|\le \kappa_B.$$
\end{ipotesi}

By Assumption \ref{assf} there
is a  $t\in (0,1)$ such that
$$ f_{\Nkp1}(x_{k}+p_k)-m_k(p_k) =\int_0^1 \left( \nabla f_{\Nkp1}(x_k+tp_k)   - \nabla f_{\Nkp1}(x_k)  \right)^T p_k  dt
 - \frac{1}{2}p_k^T B_{k+1}p_k , $$
\cite[Lemma 4.1.2]{ds}. Consequently, using  Assumptions \ref{assxk}--\ref{Bk} we have 
\begin{equation}\label{diff_fm}
|f_{\Nkp1}(x_{k}+p_k)-m_k(p_k) |\le \kT \Delta_k^2,
\end{equation}
with  $\kT=(L+\kappa_B/2)$ and  $L$ depending on the Lipschitz constants of $ \nabla \phi_i$, $1\le i\le N$.

In the next result we use the key inequality
\be \label{eqn8}
m_k(0) - m_k(p_k^C) \geq \frac 1 2 \|\nabla f_{\Nkp1}(x_k)\| \min \left\{\frac{\|\nabla f_{\Nkp1}(x_k)\|}{\beta}, \Delta_k\right\},
\ee
with $ \beta=1 +\kappa_B$, see \cite[Theorem 6.3.1]{cgt}.

\begin{lemma}
Let Assumptions \ref{assh}-- \ref{Bk} hold. Assume  $\t_k\in (0,1)$ and $\t_{k+1}$  as in (\ref{tkp1}). Then, 
  Steps 4 and 5 of Algorithm {\sc IRETR}  are well defined.
\end{lemma}
\bpr 
Let us  
prove that $\ared(\t_{k+1})-\eta \pred(\t_{k+1})$ is strictly positive if   $\Delta_{k}^{({\cal T}_k)}$ is small enough,   i.e., after a finite number $ {\cal T}_k$ of reductions of the trust-region radius. 
Let $ \theta_{k+1} $ be computed at Step 4 for some $ \Delta_k^{({\cal T}_k)}. $ 
By (\ref{defpred}) and (\ref{defared}), we have
\begin{eqnarray*}
&&\ared(\t_{k+1})-\eta \pred(\t_{k+1})\\ 
& & = \t_{k+1} (f_{\Nk}(x_k)-f_{\Nkp1}(x_{k}+p_k))+(1-\t_{k+1})(h(\Nk)-h(\Nkp1))-\\
& & \quad \eta\t_{k+1}(f_{\Nk}(x_k)-m_k(p_k))-\eta(1-\t_{k+1}) (h(\Nk)-h(\Ntilde))\\
& & =  \t_{k+1} (f_{\Nk}(x_k)-m_k(p_k))+ \t_{k+1}(m_k(p_k)-f_{\Nkp1}(x_{k}+p_k))+\\
& & \quad(1-\t_{k+1})(h(\Nk)-h(\Ntilde))+ (1-\t_{k+1})(h(\Ntilde)-h(\Nkp1))-\\
& &\quad \eta\t_{k+1}(f_{\Nk}(x_k)-m_k(p_k))-\eta(1-\t_{k+1}) (h(\Nk)-h(\Ntilde))\\
& &= (1-\eta) \left(\t_{k+1}(f_{\Nk}(x_k)-m_k(p_k)) +(1-\t_{k+1})(h(\Nk)-h(\Ntilde)) \right)+\\
& &\quad \t_{k+1}(m_k(p_k)-f_{\Nkp1}(x_{k}+p_k)) +  (1-\t_{k+1})(h(\Ntilde)-h(\Nkp1))\\
& & =(1-\eta)  \pred(\t_{k+1}) + \t_{k+1}(m_k(p_k)-f_{\Nkp1}(x_{k}+p_k)) +  \\
& & \quad(1-\t_{k+1})(h(\Ntilde)-h(\Nkp1)).
\end{eqnarray*} 
We now distinguish three cases. 

$i)$ If $h(\Nk)-h(\Ntilde)>0$ then using (\ref{pred}) we get 
\begin{eqnarray} 
 \ared(\t_{k+1})-\eta \pred(\t_{k+1}) &\ge& \eta (1-\eta) (h(\Nk)-h(\Ntilde)) \nonumber\\
 &&
  +\t_{k+1}(m_k(p_k)  -f_{\Nkp1}(x_{k}+p_k)) \nonumber\\ 
& & +   (1-\t_{k+1})(h(\Ntilde)-h(\Nkp1)) . \label{bound}
\end{eqnarray} 
The first term in the above right hand-side is  strictly  positive and uniformly bounded from below due to (\ref{riduh}).
On the other hand,  by (\ref{diff_fm}) and (\ref{new2})  
\begin{equation} \label{bound_1}
 |\t_{k+1}(m_k(p_k)-f_{\Nkp1}(x_{k}+p_k)) +   (1-\t_{k+1})(h(\Ntilde)-h(\Nkp1)) |
 \le  \kT\left(\Delta_k^{({\cal T}_k)}\right)^2+\mu \left(\Delta_k^{( {\cal T}_k)}\right)^{1+\gamma}
\end{equation} 
Therefore, for $\Delta_k^{{\cal T}_k}$   small enough we have  $ \ared(\t_{k+1})-\eta \pred(\t_{k+1}) > 0  $ and the iteration finishes.

$ii)$ 
If  $h(\Nk)-h(\Ntilde)=0$ (equivalently $\Nk=\Ntilde=N$) and $\Nkp1=N$  then  using
(\ref{defpred}) and (\ref{defared}) we have
\begin{eqnarray*} 
 \ared(\t_{k+1})-\eta \pred(\t_{k+1})&=&
 (1-\eta)\t_{k+1} (f_{N}(x_k)-m_k(p_k)) +\\
&& \quad \t_{k+1}(m_k(p_k)-f_{N}(x_{k}+p_k) ).
\end{eqnarray*}
Thus, by  (\ref{steptr}), (\ref{diff_fm}) and  (\ref{eqn8}), if  $ \Delta_k^{({\cal T}_k)}$ is small enough  we get
\begin{eqnarray}
 & & \ared(\t_{k+1})-\eta \pred(\t_{k+1})  \nonumber \\
& & \quad\ge  \tau (1-\eta)\t_{k+1} (m_k(0)-m_k(p_k^C)) +\t_{k+1}(m_k(p_k)-f_{N}(x_{k}+p_k) )\nonumber \\
& &\quad \ge \frac{1}{2}\tau (1-\eta)\t_{k+1} \|\nabla f_N(x_k)\| \Delta_k^{({\cal T}_k)}-\t_{k+1}|m_k(p_k)-f_{N}(x_{k}+p_k) | \nonumber \\
& & \quad\ge \frac{1}{2}\tau \underline{\t}(1-\eta)\|\nabla f_N(x_k)\| \Delta_k^{({\cal T}_k)}-|m_k(p_k)-f_{N}(x_{k}+p_k) | \nonumber \\
& & \quad\ge \left (\frac{1}{2}\tau  \underline{\t} (1-\eta)\|\nabla f_N(x_k)\|-\kT \Delta_k^{({\cal T}_k)}\right) \Delta_k^{({\cal T}_k)}, 
\label{set2} 
\end{eqnarray} 
and the last bound is positive   for some finite  ${\cal T}_k$. 

$iii)$
Finally,  suppose  $h(\Nk)-h(\Ntilde)=0$   (equivalently $\Nk=\Ntilde=N$) and $\Nkp1<N$   then using (\ref{defpred}) and (\ref{defared}) 
we have 
\begin{eqnarray*}
& & \ared(\t_{k+1})-\eta \pred(\t_{k+1})=  (1-\eta)\t_{k+1} (f_{N}(x_k)-m_k(p_k)) +\\
&& \hspace*{100pt}  \t_{k+1}(m_k(p_k)-f_{\Nkp1}(x_{k}+p_k))
- (1-\t_{k+1})h(\Nkp1).
\end{eqnarray*} 
Thus, by Step 3 of Algorithm 2.1, (\ref{diff_fm})  and   (\ref{eqn8}), if  $ \Delta_k^{{\cal T}_k}$ is small enough we get
\begin{eqnarray}
 & & \ared(\t_{k+1})-\eta \pred(\t_{k+1}) \ge \\
  & & \quad \ge  (1-\eta)\underline{\t}\tau  (m_k(0)-m_k(p_k^C)) - \t_{k+1}|m_k(p_k)-f_{\Nkp1}(x_{k}+p_k) |- h(\Nkp1) \nonumber \\
& & \quad \ge  \frac{1}{2}\tau\underline{\t} (1-\eta) \|\nabla f_{\Nkp1}(x_k)\|\Delta_k^{({\cal T}_k)}-\t_{k+1}|m_k(p_k)-f_{\Nkp1}(x_{k}+p_k) |- 
h(\Nkp1) \nonumber \\
& &\quad  \ge   \left( \frac{1}{2}\tau\underline{\t} (1-\eta)\|\nabla f_{\Nkp1}(x_k)\|- \kT \Delta_k^{({\cal T}_k) } - \mu\left(\Delta_k^{({\cal T}_k)}\right)^{ \gamma}\right)\Delta_k^{({\cal T}_k)}, \label{set3} 
\end{eqnarray} 
and the last bound is positive for some finite ${\cal T}_k$.
\epr

The analysis presented in the rest of this section concerns the case where Algorithm {\sc iretr}   is invoked with $\varepsilon_g=0$ and does not terminate in a finite number of steps.
Each iteration $ k-1 $ of the Algorithm ends up with the accepted iterate $ x_{k} = x_{k-1}+p_{k-1}  $ and the final sample size $ N_{k}. $ In the following statements we are going to prove that $ h(N_k) \to 0 $ and therefore the full sample is eventually reached and maintained. 
\begin{theorem}\label{limh}
Let Assumptions \ref{assh}--\ref{Bk}   hold. 
Then $ h(N_k) \to 0$.
\end{theorem}
\bpr
Inequalities  (\ref{feas}) and (\ref{pred}) imply   
\begin{equation}\label{hbound}
 h(\Nk)\le \frac{h(\Nk)-h(\Ntilde)}{1-r}\le \frac{\pred(\t_{k+1})}{\eta(1-r)}. 
 \end{equation}
We prove by contradiction that $\lim_{k \rightarrow \infty} \pred(\t_{k+1})=0$. 

 Taking into account that at termination of iteration $k$ we have $x_{k+1}=x_k+p_k$ and $\ared(\t_{k+1})\ge \eta\pred(\t_{k+1})$,  using (\ref{ared}) and (\ref{defared}) we have
\begin{eqnarray*}
 {\t_{k+2}f_{\Nkp1}(x_{k+1})} &\le & \t_{k+2}f_{\Nkp1}(x_{k+1})+\ared(\t_{k+1})-\eta\pred(\t_{k+1})\\
&=&\t_{k+1}f_{\Nk}(x_k)+ (\t_{k+2}-\t_{k+1})f_{\Nkp1}(x_{k+1})+\\
& & (1-\t_{k+1})(h(\Nk)-h(\Nkp1))-  \eta\pred(\t_{k+1})\\
&\le& \t_{k+1}f_{\Nk}(x_k)+(\t_{k+1}-\t_{k+2})\max_{x\in \Omega} \sum_{ i\in I_{\Nkp1} }|\phi_i(x)|+ \\
& &  (1-\t_{k+1})(h(\Nk)-h(\Nkp1))-\eta\pred(\t_{k+1}).
\end{eqnarray*}
Using (\ref{kappaphi})  we can rewrite the above inequality as 
\begin{eqnarray*}
   \t_{k+2}f_{\Nkp1}(x_{k+1}) & \le & \t_{k+1}f_{N_k}(x_{k}) + (\t_{k+1}-\t_{k+2})N\kappa_{\phi} \\ 
 &+ & (1-\t_{k+1})(h(\Nk) 
  -h(\Nkp1))- \eta\pred(\t_{k+1}) .
\end{eqnarray*}
Then using recurrence,   and $-(1-\t_{k+1}) h(\Nkp1)\le 0$ we get
\begin{eqnarray*}
\t_{k+2}f_{\Nkp1}(x_{k+1})&\le& \t_{k} f_{N_{k-1}}(x_{k-1})+(\t_{k}-\t_{k+2})N\kappa_{\phi}+ 
   (1-\t_k)(h(N_{k-1})-h(N_k))\\ & & +(1-\t_{k+1}) (h(N_k)-h(N_{k+1}) - 
\eta\sum_{j=k-1}^{k}{\rm Pred}_{j}(\t_{j+1}) \\
&\le&    \t_{k} f_{N_{k-1}}(x_{k-1})+(\t_{k}-\t_{k+2})N\kappa_{\phi}+\\
 & &  (1-\t_k)h(N_{k-1})+(\t_k-\t_{k+1}) h(N_k)-\eta\sum_{j=k-1}^{k}{\rm Pred}_{j}(\t_{j+1}).  
\end{eqnarray*}
Repeating this argument, using $(\t_{j}-\t_{j+1})\ge 0$ from Lemma \ref{lemmatheta} and  (\ref{boundh})  we  obtain
\begin{eqnarray*}
\t_{k+2}f_{\Nkp1}(x_{k+1})&\le& \t_1 f_{N_0}(x_0)+(\t_{1}-\t_{k+2})N\kappa_{\phi}+ (1-\t_1)h(N_0) +\\
& &  \sum_{j=1}^k (\t_j-\t_{j+1}) h(N_j) -\eta\sum_{j=0}^{k}{\rm Pred}_{j}(\t_{j+1}) \\
  &\le&    \t_1 f_{N_0}(x_0)+ (1-\underline{\t})N\kappa_{\phi}+ (1-\underline{\t})\bar h+\sum_{j=1}^k (\t_j-\t_{j+1}) \bar h  -\\
	& & \eta\sum_{j=0}^{k}{\rm Pred}_{j}(\t_{j+1})\\
		  &\le&   \t_1    f_{N_0}(x_0) + (1-\underline{\t})N\kappa_{\phi}+ (1-\underline{\t})\bar h +  (\t_1-\t_{k+1}) \bar h  -\\
	& & \eta\sum_{j=0}^{k}{\rm Pred}_{j}(\t_{j+1})\\
	  &\le&  \t_1 f_{N_0}(x_0)+ (1-\underline{\t}) N\kappa_{\phi} +   2 (1-\underline{\t})\bar h   -\eta\sum_{j=0}^{k}{\rm Pred}_{j}(\t_{j+1}).
\end{eqnarray*}
By  (\ref{newassf})    and (\ref{boundh}) we have 
\begin{eqnarray*}
\t_{k+2}f_{\Nkp1}(x_{k+1}) &=& \t_{k+2}f_{N}(x_{k+1})+\t_{k+2}(f_{\Nkp1}(x_{k+1})-f_{N}(x_{k+1}) )\\
&\ge& \t_{k+2}f_{N}(x_{k+1})-\t_{k+2}|(f_{\Nkp1}(x_{k+1})-f_{N}(x_{k+1}) )| \\
&\ge& \t_{k+2}f_{N}(x_{k+1})-\sigma \bar h,
\end{eqnarray*}
and  therefore 
\be\label{sommapred}
 \t_{k+2}f_{N}(x_{k+1}) \le \xi-\eta\sum_{j=0}^{k}{\rm Pred}_{j}(\t_{j+1}),
\ee
where 
\begin{equation}\label{xi}
\xi= \t_1 f_{N_0}(x_0)+ (1-\underline{\t})N\kappa_{\phi}+2 (1-\underline{\t}+\sigma)\bar h ,
\end{equation} 
is independent of $k$. 

  Noting that ${\rm Pred}_{j}(\t_{j+1})\ge 0$, we can conclude that if  ${\rm Pred}_{j}(\t_{j+1})$ is not tending to zero, then $\sum_{j=0}^{\infty}{\rm Pred}_{j}(\t_{j+1})$ is diverging  and this implies that 
$f_{N}$ is unbounded below in $\Omega$. This contradicts the compactness of $\Omega$.
\epr

\begin{corollary}\label{fullprecision}
Let Assumptions \ref{assh}--\ref{Bk}   hold.
Then  $N_k=N$ for all $k$ sufficiently large.
\end{corollary}
\bpr
By Theorem \ref{limh} and Assumption \ref{assh}, it follows $h(\Nk)<h(N-1)$ for all $k$ sufficiently large. This implies $\Nk=N$.
\epr

\begin{corollary}\label{first-order}
Let Assumptions \ref{assh}--\ref{Bk}   hold. 
Then, for $k$ sufficiently large, the iterations are generated by a (standard) trust-region scheme on $f_{N}$ and 
\begin{description}
\item{i)} $\liminf_{k\rightarrow \infty} \|\nabla f_N(x_{k})\|=0$.
\item{ii)} $\lim_{k\rightarrow \infty} \|\nabla f_N(x_{k})\|=0$, provided that $f_N$ is Lipschitz continuous in $\Omega$.
\end{description}
\end{corollary}
\bpr
By Corollary \ref{fullprecision} we know that  at termination of iteration $k-1$  we have  $N_k=N$ for all $k$ sufficiently large. Thus eventually,
$x_{k+1}=x_k+p_k$ with $p_k$ satisfying (\ref{ared}) which now takes the form of the standard acceptance rule of the trial point 
in trust-region methods, i.e, 
$$
\frac{f_{N}(x_{k+1})-f_N(x_k)}{m_k(0)-m_k(p_k)}\ge \eta.
$$
As a consequence,   Theorem 4.6 in \cite{nocedal} yields item $i)$. 
Item  $ii)$ 
is guaranteed by 
\cite[Theorem 4.7]{nocedal}.
\epr


\section{On the realization of the algorithm} \label{Sec3}
The realization of Algorithm {\sc iretr} raises many issues and in this section we discuss two important aspects:  the form of the model 
used and related properties, and a computationally convenient  adaptation of the rule for choosing $\Nkp1$ eventually.
We will further address implementation issues in Section \ref{exps}.

Various models of the form (\ref{model}) can be built. One possibility is the linear model 
$$ m_k(p) = f_{\Nkp1}(x_k) + \nabla f_{\Nkp1}(x_k)^T p,$$
which gives rise to a gradient method and step   $p_k$ 
$$ p_k = - \Delta_k \frac{\nabla f_{\Nkp1}(x_k)}{\|\nabla f_{\Nkp1}(x_k)\|}. $$
Namely, Algorithm {\sc iretr} becomes  a  subsampled  gradient method with variable stepsize  determined accordingly to the trust-region strategy. 

Another possibility is to use quadratic models
of the form 
$$ m_k(p) =  f_{\Nkp1}(x_k) + \nabla f_{\Nkp1}(x_k)^T p + \frac{1}{2} p^T B_{k+1} p, $$ 
and fully exploit the advantages of the trust-region framework. 
If all functions $ \phi_i $ are  twice continuously differentiable one can build the quadratic model 
%
%
$$ m_k(p) = f_{\Nkp1}(x_k) + \nabla f_{\Nkp1}(x_k)^T p + \frac{1}{2} p^T \nabla^2 f_{\Dkp1}(x_k) p, $$
with $1\le \Dkp1\le \Nkp1$ and $I_{\Dkp1}\subseteq I_{\Nkp1}$. In fact, 
the Hessian matrix $\nabla^2 f_{\Nkp1}(x)$ is approximated via subsampling by
\begin{equation}\label{Bk_sampling}
B_{k+1}=\frac{1}{\Dkp1}\sum_{i\in I_{\Dkp1}}\nabla^2 \phi_i(x_k).
\end{equation} 
The cardinality of $ I_{\Dkp1} $ now controls the precision of Hessian approximation and allows for trade-off 
between precision and computational costs. 
This particular form of Hessian approximation will be analysed in details for strongly convex functions in the next section.

The use of  quadratic models  is crucial for the computation of  
$(\varepsilon_g, \varepsilon_H)$-approximate second order critical point  of nonconvex problems  (\ref{minf}), i.e., a point $x$ such that
\begin{equation}\label{approx_opt}
\|\nabla f_{N}(x)\|\le \varepsilon_g, \quad \lambda_{\min}(\nabla^2 f_{\Dkp1}(x))\ge -\varepsilon_H,
\end{equation}
Supposing that full precision is reached, $\Nk=N$, the   trust-region problem (\ref{tr_pb}) has to be solved approximately finding 
$p_k$ such that
\begin{eqnarray}\label{eig}
m_k(p_k)\le m_k(p_k^C)  \, \mbox{ and } \,   m_k(p_k)\le m_k(p_k^E)  \  \  \mbox {if } \   \lambda_{\min}( \nabla^2 f_{\Dkp1} (x_k))<0,
\end{eqnarray} 
where $p_k^C$ is the  Cauchy point (\ref{cauchy}) and $p_k^E$ is a negative curvature direction such that  $(p_k^E)^T \nabla^2 f_{\Dkp1} (x_k) p_k^E\le\upsilon  
\lambda_{\min}( \nabla^2 f_{\Dkp1} (x_k) )\|p_k^E\|^2$ for some $\upsilon\in (0,1]$, \cite[\S 6.6]{cgt}. 

We refer to \cite[Theorem 1]{XuRoosMaho17} for results on the computation of approximated  second-order optimal solutions using 
trust-region methods with full function  and gradient  and subsampled Hessian.
%

  Let us now address the choice   of  the stopping criterion in Algorithm {\sc iretr}. 
Notice that the Algorithm may stop even if  full precision  at iteration $k$ is not achieved (i.e. $\Nkp1<N$), provided that 
  $N_k=N$. This choice  is supported by observing that 
suitable sample sizes provide an accurate approximation  $\nabla f_{\Nkp1}(x_k)$ to $\nabla f_N(x_k)$.
In fact, by \cite[Theorem 6.2]{bgmt}  $\nabla f_{\Nkp1}(x_k)$ is sufficiently accurate with  fixed probability at least  $1-p_g$, i.e., 
$$
Pr(\| \nabla f_N(x_k)-\nabla f_{\Nkp1}(x_k) \|\le \chi_g )\ge 1-p_g \ \ \mbox{with} \ \   \chi_g\in (0,1), \,  \  p_g\in (0,1),
$$
if the cardinality $\Nkp1$  satisfies
\begin{equation}\label{SK1}
\Nkp1\ge \min\left \{ N,\left\lceil  \frac{2 }{\chi_g } \left( \frac{ V_{g}}{\chi_g }+\frac{2\zeta(x_k)}{3}\right) 
  \,\log\left(\frac{n+1}{ p_g}\right)\right\rceil\right \} ,
\end{equation}
with  $E(\| \nabla \phi_i(x_k)-\nabla f_N(x_k)\|^2)\le V_g$ and $\max_{i\in\{1,...,N\}} |  \nabla  \phi_i(x) | \le \zeta(x)$, 
and $I_{\Nkp1}$ is sampled uniformly in $\{1, 2, \ldots, N\}$.


We conclude this section observing that, in the current form of the algorithm, at each iteration an attempt is made to use $\Nkp1<N$ (see Step 2). By Corollary 
\ref{fullprecision} we know that, for $k$ sufficiently large, such a value will be rejected and this fact implies  
useless repetitions of Steps 2--5. To overcome this drawback, we replace (\ref{new2}) with
\begin{eqnarray}
 h(\Nkp1)-h(\Ntilde) &\le& \mu \Delta_k^{1+\gamma} \hspace*{101 pt}  \mbox{ if }  \ \ \Nk\neq N  \label{new31}\\
 h(\Nkp1)-h(N) &\le& \min\{ \mu \Delta_k^{1+\gamma}, \, \|\nabla f_N(x_k) \| \} \quad \mbox{ if } \ \ \Nk= N  \label{new32} 
\end{eqnarray}
Then, the following result holds.
\begin{corollary} \label{C2.7}
Suppose (\ref{new31}) and (\ref{new32}) hold. 
For $k$ sufficiently large,  the use of  sets $I_{\Nkp1}$ of cardinality smaller than $N$ is  not 
attempted.
\end{corollary}
\bpr
By Corollary \ref{fullprecision}  and Corollary \ref{first-order}, 
we know that $N_k=N$ for all $k$ sufficiently large and  $\|\nabla f_N(x_{k})\|$ tends to zero. Thus,
letting $k_*$ be the iteration index such that  $\|\nabla f_N(x_{k})\|<h(N-1)$, $\forall k\ge k_*$, it follows
$\Nkp1=N$, $\forall k\ge k_*$.
\epr

\section{Strongly convex problems}
In this section we assume that $f_N$ is strongly convex  with strongly convex functions $\phi_i$, $1\le i\le N$,
and  analyze the local behaviour  of {\sc iretr} method when 
full precision for the function and the gradient has been reached and a quadratic model  of  the following form is used:
$$ m_k(p) = f_{N}(x_k) + \nabla f_{N}(x_k)^T p + \frac{1}{2} p^T \nabla^2 f_{\Dkp1}(x_k) p, $$
with $1\le \Dkp1\le N$, $I_{\Dkp1}\subseteq I_N$. 
Thus, we are focusing on the local behaviour of the trust-region method  
employing second order models with exact function and gradient and subsampled Hessian.
Such a method has been investigated in \cite{XuRoosMaho17} with respect to iteration complexity
but not with respect to local convergence.

The additional assumptions used in this section are stated below. 
\begin{ipotesi} \label{A1} The functions $ \phi_i, \; i=1,\ldots, N $, are twice continuously differentiable and strongly convex in $\IR^n$, 
\be \label{eqn1}
 \lambda_1  I \preceq \nabla^2 \phi_i(x) \preceq \lambda_n  I, \; \mbox{ with  } 0 < \lambda_1 < \lambda_n,
 \ee
where, given two matrices $A$ and $B$,  $A \preceq B$ means that $B-A$ is positive semidefinite.
\end{ipotesi}

Trivially, $f_N$ is strongly convex and admits an unique minimizer $x^*$. Moreover, 
 $B_{k+1}$ is as in \eqref{Bk_sampling}, both 
 $\lambda_{\min}(B_{k+1})\ge \lambda_1$ and $\lambda_{\max}\le \lambda_n$  hold and Corollary \ref{first-order} implies 
$\lim_{k\rightarrow \infty}  x_k=  x^*$.

The following theorem analyzes the behaviour of $\{x_k\}$ denoting  
\begin{equation}\label{errD}
e(\Dkp1) =  \|\nabla^2 f_N(x_k) - \nabla^2 f_{\Dkp1}(x_k)\|,
\end{equation}
the error between  $\nabla^2 f_N(x_k)$ and  $\nabla^2 f_{\Dkp1}(x_k) $. 
We also invoke the assumption below.
\begin{ipotesi} \label{A2} The Hessian $\nabla^2 f_N $ is Lipschitz continuous on 
$ {\cal B}_{\delta}(x^*):=\{x \in \mathbb{R}^n: \|x-x^*\| \leq \delta\} $ with Lipschitz constant $2 L_H $.
\end{ipotesi}

\begin{theorem} \label{T2}
Suppose that  Assumptions  \ref{assh}, \ref{assxk}, \ref{A1},  \ref{A2} hold.
Let $ \{x_k\}$ be generated by Algorithm {\sc iretr}, $\varepsilon_g$ as in (\ref{epsfo}), $\beta$ as in (\ref{eqn8}),
$\eta$ as in  the Algorithm {\sc iretr} and $B_{k+1}$ given by \eqref{Bk_sampling}.

\begin{description}
\item{i)}  Let  $ \epsilon\in (0,1) $   and $ \Dkp1 $ such that
\be \label{eqn9} 
\displaystyle \frac{1}{{\tau}\min\left\{ \frac{\lambda_1^2}{4\beta}, \frac {\lambda_1}{ 2} \right\}} \left(2 L_H\epsilon+ e(\Dkp1) \right) \leq 1 - \eta.
\ee 
Then,  if $k$ is sufficiently large,  $p_k$ is accepted in the first pass in Step 5 and  $ {\cal T}_k = 0. $
 \item{ii)}  There exist sufficiently small  $ \delta > 0 $ and sufficiently large $ D $ such that, for all $ x_k \in {\cal B}_{\delta}(x^*)$
and $\Dkp1=D, $  the error $\|x_k-x^*\|$ reduces linearly, i.e., 
$\|x_{k+1}-x^*\|< \tilde{\tau} \|x_k-x^*\|$  for some $\tilde{\tau} \in (0,1)$. 
\end{description}
\end{theorem}
\bpr
$i)$
 Let us consider $ k $ sufficiently large such that $ N_{k+1} = N $ at termination of iteration $k$. Lemma 6.5.1 in \cite{cgt}  gives 
\begin{equation}\label{step}
\|p_k\| \leq \frac{2}{\lambda_1} \|\nabla f_N(x_k)\|.
\end{equation}
 Let us consider the step $p_k$ returned by iteration $k$.
  Combining  \eqref{step} with (\ref{eqn8}) and \eqref{steptr} we obtain 
\be \label{eqn81}
m_k(0) - m_k(p_k) \geq \frac 1 2  \omega \|p_k\|^2,
\ee
with $\omega= \tau \min\{\frac{\lambda_1}{2\beta},1\} \frac{\lambda_1}{2}$. 

At Step 5 of the Algorithm, (\ref{ared}) has the form $f_N(x_k)-f_N(x_k+p_k) \ge \eta (m_k(0) - m_k(p_k))$.
By Assumption \ref{A2} and (\ref{errD}), it follows
\begin{eqnarray*}
  \left| \frac{f_N(x_k)-f_N(x_k+p_k) }{m_k(0) - m_k(p_k)}-1  \right|
  &=&\frac{|f_N(x_k+p_k)-m_k(p_k)|}{ m_k(0) - m_k(p_k) }\\
&  \le& \frac{|\frac{1}{2}p_k^T (\nabla^2 f_N(x_k+tp_k)-\nabla^2 f_{\Dkp1}(x_k))p_k|}{\frac 1 2\omega \|p_k\|^2}\\
& \le& \frac{1}{2}\|p_k\|^2\Big(\frac{\|\nabla^2 f_N(x_k+tp_k)-\nabla^2 f_{N}(x_k)\|}{\frac 1 2 \omega \|p_k\|^2}\\ 
&& +
\frac{\|\nabla^2 f_N(x_k)-\nabla^2 f_{\Dkp1}(x_k)\|}{\frac 1 2 \omega \|p_k\|^2}\Big)\\
&   \le& \frac{2 L_H\|p_k\|+ e(\Dkp1)}{\omega},
\end{eqnarray*}
where $t$ is some scalar in $t\in (0,\, 1)$ \cite[Theorem 3.1.2]{cgt}.
Now, given $\epsilon\in (0,1)$ and $D_{k+1}$ as in (\ref{eqn9}), (\ref{step})  and Corollary \ref{first-order} imply $\|p_k\| \le \epsilon$ for $k$ large enough, say $k\ge \bar k$,  and (\ref{eqn9}) implies the acceptance of the step. Then, $\Delta_k$ is not reduced and $\Delta_k\ge \Delta_{\bar k}$  for any $k\ge \bar k$.

$ii)$  
Using  \eqref{step}, Corollary \ref{first-order} and item $i)$ we can conclude 
that the trust-region bound becomes inactive for $k$ sufficiently large, 
i.e.,  the step 
$$ p_k = -( \nabla^2 f_{\Dkp1}(x_k))^{-1} \nabla f_{N}(x_k) ,$$
is accepted eventually. Consequently, using multivariate calculus results \cite[Lemma 4.1.12]{ds} and Assumption \ref{A1}
\begin{eqnarray}
 \|x_{k+1}-x^*\| &=&\|x_k - ( \nabla^2 f_{\Dkp1}(x_k))^{-1} \nabla f_{N}(x_k)  -x^*\| \nonumber \\ 
& =&  \|( \nabla^2 f_{\Dkp1}(x_k))^{-1}  ( \nabla f_N(x^*)- \nabla f_{N}(x_k)  - \nabla^2 f_{\Dkp1} (x_k)(x^*-x_k)\|\nonumber \\ 
& \le&  \|( \nabla^2 f_{\Dkp1}(x_k))^{-1}\| \, \left(\| \nabla f_N(x^*)- \nabla f_{N}(x_k)  - \nabla^2 f_{N} (x_k)(x^*-x_k)\|\right. \nonumber \\ 
& &  \left.+\|( \nabla^2 f_{N} (x_k)-\nabla^2 f_{\Dkp1}(x_k))(x^*-x_k)\|\right)\nonumber \\
&  \le& \frac{1}{\lambda_1} \|x_k-x^*\|\left( L_H\|x_k-x^*\|+ e(\Dkp1)\right) \label{bound1} 
\end{eqnarray}
Thus, the claim follows if $\delta$ and $\Dkp1=D$ are such that 
$\tilde{\tau} := \frac{L_H\delta+e(D)}{\lambda_1}<1$ and $D$ satisfies (\ref{eqn9}).

\epr

 Item $ii)$ above may require a rather large value $\Dkp1=D$
which is adverse for practical computation.
A  more stringent condition on  $D_{k+1} $ of the form 
$e(D_{k+1})=O(\|\nabla f_N(x_k)\| )$  yields quadratic convergence but again  such $ D_{k+1} $ might be very close to $ N $.
We next investigate  on the more realistic situation where  the Hessian accuracy requirement in  \eqref{eqn9} is guaranteed only with high-probability and provide a linear convergence result in expectation.

Let us now  suppose that, given an accuracy requirement $ \chi_H>0$,
the probability of $\|\nabla^2 f_N(x_k)-\nabla^2 f_{\Dkp1}(x_k)\|$ being smaller than $ \chi_H $ is larger than $ 1 - p_H$:
\be \label{eqn12} 
P(\|\nabla^2 f_N(x_k) - \nabla^2 f_{\Dkp1}(x_k)\| \leq \chi_H) \geq 1 - p_H,
\ee 
for $p_H \in (0,1)$.
If the subsample $I_{\Dkp1}$ is chosen randomly and uniformly, then  
the   lower bound on the  sample size ensuring \eqref{eqn12} takes the form
\be \label{minbound} 
 D_{k+1} \geq  \min \left\{ N,\left\lceil\frac{ 2}{\chi_H}\left( \frac{\lambda_n^2}{ \chi_H}+\frac{\lambda_n}{3}\right)\log\left(\frac{2n}{p_H}\right)\right\rceil  \right\}.  
\ee
The above bound  is derived  in \cite[Lemma 3.1]{bkkj}   and a similar bound is given in  \cite[Lemma 4]{bgm}.

We now provide a linear convergence result in expectation; the   step $p_k$ taken is the global minimizer of  (\ref{tr_pb}), i.e., 
$$
(\nabla^2 f_{\Dkp1}(x_k)+\nu_k I)p_k=-\nabla f_N(x_k),
$$
for some $\nu_k\ge 0$, see \cite[Theorem 7.2.1]{cgt}.

\begin{theorem} \label{Texpectation} 
Suppose that  Assumptions  \ref{assh}, \ref{assxk}, \ref{A1}, \ref{A2}  hold.  Let $ \{x_k\}$ be generated by Algorithm {\sc iretr}  invoked with $\varepsilon_g=0$ in (\ref{epsfo}),  $B_{k+1}$ as in  \eqref{Bk_sampling}
and $p_k$ being the global minimizer of (\ref{tr_pb}).
If (\ref{eqn12}) holds and there exists a  $\nu^*\in (0,1)$ such that for all $k$
\begin{equation}\label{glob}
\frac{\nu_k}{\lambda_1+\nu_k}\le \nu^* ,
\end{equation}
 then there exist $\delta$, $ \chi_H   $, $ p_H  $ sufficiently small  such that  
\be \label{mslinear} 
E(\|x_{k+1} - x^*\|) \leq \bar \tau E(\|x_{k} - x^*\|), 
\ee 
for all $ k $ large enough and some $ \bar\tau \in (0,1). $  
\end{theorem}
\bpr
 Take $ \delta \in (0,1) $, $ \chi_H >0 $, $p_H\in (0,1)$ such that 
\begin{eqnarray}
& &  \rho=  \frac{L_H \, \delta}{\lambda_1 }  + \frac{\chi_H}{\lambda_1 }+\nu^* \le \bar \tau ,\label{rho1}  \\
& &  
p_H \leq \frac{(\bar \tau   - \rho)}{1 + \frac{2\lambda_1}{\lambda_n} } .\label{rho2} 
\end{eqnarray}
for some $ \bar \tau \in (0,1). $ 
Let   $ k$ large enough such that  $x_k \in   {\cal B}_{\delta}(x^*)$.

Denote by $ A_k $ the event 
\be \label{event_pa}
  \|\nabla^2 f_{\Dkp1}(x_k) - \nabla^2 f_N(x_k)\| \leq \chi_H. 
  \ee
Then $ P(A_k) \geq 1 - p_H $ and $ P(\bar{A}_k) < p_H,$  where $  \bar{A}_k $ denotes the event  $A_k$ does not occur.
If $ A_k $ happens then  using multivariate calculus results \cite[Lemma 4.1.12]{ds},  Assumption \ref{A1},
(\ref{glob}) and  (\ref{rho1})
\begin{eqnarray}
 \|x_{k+1}-x^*\| &=&\|x_k - ( \nabla^2 f_{\Dkp1}(x_k)+\nu_k I)^{-1} \nabla f_{N}(x_k)  -x^*\| \nonumber \\ 
& =&  \|( \nabla^2 f_{\Dkp1}(x_k)+\nu_k I)^{-1}  ( \nabla f_N(x^*)- \nabla f_{N}(x_k)  - (\nabla^2 f_{\Dkp1}(x_k)+\nu_k I)(x^*-x_k)\|\nonumber \\ 
& \le&  \|( \nabla^2 f_{\Dkp1}(x_k)+\nu_k I)^{-1}\| \, \left(\| \nabla f_N(x^*)- \nabla f_{N}(x_k)  - \nabla^2 f_{N} (x_k)(x^*-x_k)\|\right. \nonumber \\ 
& &  \left.+\|( \nabla^2 f_{N} (x_k)-\nabla^2 f_{\Dkp1}(x_k))(x^*-x_k)+ \nu_k(x^*-x_k)\|\right)\nonumber \\
&  \le& \frac{1}{\lambda_1+\nu_k}\left( L_H\|x_k-x^*\|+ e(\Dkp1)+\nu_k\right) \|x_k-x^*\| \nonumber \\ 
&\le&  \left(\frac{L_H \, \delta}{\lambda_1 }  + \frac{\chi_H}{\lambda_1 }+\nu^*\right)  \|x_k - x^*\| \nonumber \\
& =&  \rho \|x_k - x^*\| \label{eqn19}
\end{eqnarray}
Otherwise, if $ \bar{A}_k $ is realized then by (\ref{step}) we have 
$$ \|x_{k+1} - x^*\| \leq \left(1 + \frac{2\lambda_1}{\lambda_n}\right) \|x_k - x^*\|. $$
Therefore, 
\begin{eqnarray*}
E(\|x_{k+1} - x^*\|) & = & P(A_k) E(\|x_{k+1} - x^*\| |A_k) + P(\bar{A}_k) E(\|x_{k+1} - x^*\| |\bar{A}_k) \\
& \leq & \rho E(\|x_{k} - x^*\|) + p_H \left(1 + \frac{2\lambda_1}{\lambda_n}\right) E(\|x_{k} - x^*\|) \\
& \leq & \bar\tau E(\|x_{k} - x^*\|),
\end{eqnarray*}
where we have used (\ref{rho2}) and $p(A_k)\le 1$. 
\epr

\section{Worst-case iteration and evaluation complexity to first-order critical points }
In this section we provide an upper bound on the number of iterations and function-evaluations needed 
to find an $\varepsilon_g$-accurate first-order optimality point (\ref{epsfo}).
The number of function-evaluations is intended as  the number of
evaluations of functions of the form $f_M$, for some $M\le N$.  We recall that a standard trust-region approach  shows   $ {\cal{O}}(\varepsilon_g^{-2}) $
worst-case iteration and full function  complexity for first-order optimality \cite{gy}.

 Recalling that $h(\Nk)-h(\Ntilde)=0$  is equivalent to $\Nk=\Ntilde=N$, 
consider the following partition of iteration  indices $k$: 
\begin{itemize}
\item ${\cal I}_1=\{ k\ge 0 \mbox{  s.t.  }   h(\Nk)-h(\Ntilde)>0\}$,
\item ${\cal I}_2=\{ k\ge 0 \mbox{  s.t.  } h(\Nk)=h(\Ntilde)=0, \Nkp1=N \mbox{ and } \|\nabla f_N(x_k)\|>\varepsilon_g \}$,
\item ${\cal I}_3=\{ k\ge 0 \mbox{  s.t.  } h(\Nk)=h(\Ntilde)=0, \, \Nkp1<N\mbox{ and } \|\nabla f_{\Nkp1}(x_k)\|>\varepsilon_g \}$.
\end{itemize}

The value of $\Nkp1$ may change within  iteration $k$ before acceptance of the iterate; above  $\Nkp1$ is the 
value at the end of iteration $k$, i.e., the value used for building the accepted iterate $x_{k+1}$.

Our analysis is carried out fixing $\gamma=1$ in Algorithm {\sc iretr} and the  
first result provides a lower bound on the trust-region radius  at termination of iteration $k$.
\begin{lemma} \label{lemma_Delta}
Let Assumptions \ref{assh}--\ref{Bk}   hold.  Suppose  furthermore   $\gamma=1$ in  Algorithm {\sc iretr}.
Then, 
\begin{description}
\item{i)} for any $k \in {\cal I}_1$
$$
\Delta_k \ge \min \left\{ \zeta_1 \sqrt{\frac{\eta (1-\eta)}{\kT +\mu} (1-r) \underline{h}},\, \Delta_0 \right\} ,
$$
\item{ii)} for any $k \in {\cal I}_2\cup {\cal I}_3$,
\begin{equation}\label{bound_delta}
\Delta_k \ge \min \left\{ \zeta_1 \sqrt{\frac{\z \underline{ h } \z }{\mu}} ,\,  \zeta_1 \Gamma \varepsilon_g ,\, \Delta_0  \right\},
\end{equation}
\end{description}
for some positive $\Gamma$ and $\mu$  as in the Algorithm.
\end{lemma}
\bpr
The initial $\Delta_k$ may be reduced in Steps 3 and 5 of the Algorithm. 
Step 3 is performed only if $k \in  {\cal I}_3$.

Let us consider case $i)$. 
Since $\gamma=1$ equation (\ref{bound_1}) becomes
 $$|\t_{k+1}(m_k(p_k)-f_{\Nkp1}(x_{k}+p_k)) +   (1-\t_{k+1})(h(\Ntilde)-h(\Nkp1)) |
 \le  (\kT+\mu) \Delta_k^2.$$
From \eqref{bound}, inequality \eqref{ared} is satisfied whenever
$$
\Delta_k \le \sqrt{\frac{\eta (1-\eta)}{\kT +\mu} (h(\Nk)-h(\Ntilde))}\,  .
$$
Thus, using \eqref{feas},  if 
$$
\Delta_k \le \sqrt{\frac{\eta (1-\eta)}{\kT +\mu} (1-r) \underline{h}} ,
$$
 then  \eqref{ared} holds and the claim  $i)$ follows from  the rule for decreasing  $\Delta_k$ in Step 5 of Algorithm {\sc iretr}.

Let us consider case $ii)$.
 Concerning Step 3, it is performed as long as $N_{k+1} <N$. Then,  
\eqref{new2} ensures that at termination of the loop in Steps  2--3  
$$\Delta_k\ge \zeta_1
\sqrt{\frac{\z \underline{ h } \z }{\mu}} \  .
 $$

Concerning  Step 5,  first suppose $k\in {\cal I}_2$ and  $\Delta_k\le  \varepsilon_g/\beta $ with $\beta$ as in (\ref{eqn8}). 
Using (\ref{set2})  we can conclude that if 
$$
\Delta_k \le \frac{\tau  \underline{\t} (1-\eta)}{2 \kT}\, \varepsilon_g,
$$
then  \eqref{ared} is satisfied. 

Suppose now  $k\in {\cal I}_3$ and $\Delta_k\le  \varepsilon_g/\beta$.  Using $\gamma=1$,  equation (\ref{set3}) becomes
$$
\ared(\t_{k+1})-\eta \pred(\t_{k+1}) 
 \ge   \left( \frac{1}{2}\tau\underline{\t} (1-\eta)\|\nabla f_{\Nkp1}(x_k)\|-(\kT+\mu)\Delta_k \right)\Delta_k,
$$
and if 
$$
\Delta_k\le   \frac{\tau\underline{\t} (1-\eta)}{2(\kT +\mu)} \, \varepsilon_g, 
$$
then  \eqref{ared} is satisfied. 

The upper bound on $\Delta_k$ for $k\in {\cal I}_3$ is sharper than the one obtained for $k\in {\cal I}_2$.
Then, due to the rule used to decrease $\Delta_k$ in Step 5, we can conclude that, at iteration
$k\in   {\cal I}_2\cup {\cal I}_3$,   condition \eqref{ared} is satisfied if 
\begin{equation}\label{Gamma}
\Delta_k>\zeta_1 \min\left\{\frac{1}{\beta},  \frac{\tau  \underline{\t} (1-\eta)}{2(\kT +\mu)}\right\}\varepsilon_g
\eqdef \zeta_1 \Gamma \varepsilon_g ,
\end{equation}
and  the claim follows.  \epr

\begin{theorem}\label{th_complexity}
Let Assumptions   \ref{assh}--\ref{Bk} hold.  Suppose  furthermore   $\gamma=1$ in  Algorithm {\sc iretr}  and 
let $f_{low}$ the lower bound of $f_N$ in $\Omega$. Then,
\begin{description}
\item{i)} the cardinality $|{\cal{I}}_1|$ satisfies
$$
{ |{\cal{I}}_1|\le \left\lceil \nu_1 \underline{h}^{-1}  \right\rceil,} 
$$
{with $\nu_1=\frac{ \xi-\underline{\t}f_{low}}{\eta^2(1-r)}$, $\xi$ as in (\ref{xi}), $\underline{\t}$ as in Lemma \ref{lemmatheta}, $\eta$ and $r$ as in the Algorithm {\sc iretr}}.

\item{ii)} the cardinality $|{\cal{I}}_2|+ |{\cal{I}}_3|$ satisfies 
$$
|{\cal{I}}_2|+ |{\cal{I}}_3|\le \left \{\begin{array}{ll} 
\left\lceil    \nu_2  \,  \varepsilon_g^{-2}  \right\rceil & \mbox{   if  }  \  \Gamma \varepsilon_g  \le \min \left\{  \displaystyle \sqrt{\frac{\z \underline{h}\z }{\mu}}, \frac{\Delta_0}{\zeta_1}\right\} ,\\
\nu_3   \underline{h}^{-\frac 1 2}   \varepsilon_g^{-1}  &    \mbox{   if  }  \    \displaystyle \sqrt{\frac{\z \underline{h}\z }{\mu}} \le \min \left\{  \Gamma \varepsilon_g, \frac{\Delta_0}{\zeta_1},\right\} 
\end{array} \right. 
$$
with positive 
$
\nu_2=\frac{2}{\eta\Gamma} \left(  f_{N_0}(x_0)-f_{low}+(\sigma\eta+1-\underline \t)  \frac{ \xi-\underline{\t}f_{low}}{\eta^2(1-r) }\right)
$,  $ \nu_3=\nu_2\Gamma \sqrt{\mu}$.

\end{description}
\end{theorem}
\bpr
Let us denote with $\bar k$ the last iterate of Algorithm {\sc iretr} and note that  $N_{\bar k}=N$ by definition of the algorithm.
From (\ref{sommapred}) it follows 
$$
 \sum_{k=0}^{\bar k -1}{\rm Pred}_{k}(\t_{k+1}) \le \frac{\xi-\t_{\bar k+1}f_{N}(x_{\bar k})}{\eta} \le \frac{ \xi-\underline{\t}f_{low}}{\eta},
\;\;\;\;\forall k\ge 0,
$$
and consequently  (\ref{hbound}) yields 
{
\begin{equation} \label{sum_h}
\sum_{k=0}^{\bar k-1} h(N_k) \le \frac{\xi-\underline{\t}f_{low}}{\eta^2(1-r)} .
\end{equation} 
}
Then the number of indices $k$ such that $h(N_k)>\underline{h}$ is bounded above by
$$\frac{ \xi-\underline{\t}f_{low}}{\underline{h}\eta^2(1-r)}, $$
and $i)$ follows.

Let us consider the case $k \in {\cal I}_2\cup {\cal I}_3$. Note that by  (\ref{defared}), (\ref{ared}), (\ref{defpred}), (\ref{newassf}) and \eqref{steptr}, we have
\begin{eqnarray*}
Ared_k(\t_{k+1}) &=&\t_{k+1}(f_{N}(x_k)-f_{\Nkp1}(x_{k+1}))-(1-\t_{k+1}) h(\Nkp1)\\
&\ge&\eta\t_{k+1}(f_{N}(x_k)-f_{\Nkp1}(x_k)+m_k(0)-m_k(p_k))\\
&\ge &-\sigma\eta\t_{k+1} h(\Nkp1) +\eta\t_{k+1} (m_k(0)-m_k(p_k))\\
 &\ge&  -\sigma\eta\t_{k+1} h(\Nkp1)+\tau \eta\t_{k+1} (m_k(0)-m_k(p_k^C))
\end{eqnarray*}
Then, by using \eqref{bound_delta} and \eqref{eqn8} it follows 
\be\label{ared1}
 f_{N}(x_k)-f_{\Nkp1}(x_{k+1}) + \sigma \eta  h(\Nkp1) \ge  \frac{ \tau \eta}{2}    
\min \left\{\zeta_1\Gamma \varepsilon_g, \zeta_1 \sqrt{\frac{\z \underline{h} \z }{\mu}}, \Delta_0 \right\} \varepsilon_g.
\ee
 Moreover, note that
due to the definition of 
$Ared_k(\t_{k+1})$ and inequalities 
\eqref{pred} and \eqref{ared}, 
 the following inequality holds at termination of each iteration $k\ge 0$:
\begin{equation}\label{fbound}
\frac{Ared_k(\t_{k+1})}{\t_{k+1}}=f_{N_k}(x_k)-f_{\Nkp1}(x_{k+1})+\frac{1-\t_{k+1}}{\t_{k+1}} (h(N_k)-h(\Nkp1))\ge 0
\end{equation}
Then, since  $\frac{Ared_k(\t_{k+1})}{\t_{k+1}}$ is positive,
$$
\sum_{k\in {\cal I}_2\cup {\cal I}_3}\frac{Ared_k(\t_{k+1})}{\t_{k+1}}\le \sum_{k=0}^{\bar k-1}\frac{ Ared_k(\t_{k+1})}{\t_{k+1}},
$$
and this implies
\begin{eqnarray*}
 \sum_{k\in {\cal I}_2\cup {\cal I}_3}\left( f_{N}(x_k)-f_{\Nkp1}(x_{k+1}) \right) &\le& \sum_{k=0}^{\bar k-1} \left( f_{N_k}(x_k)-f_{\Nkp1}(x_{k+1})\right )  \\ 
+ \sum_{k=0}^{\bar k-1} \frac{1-\t_{k+1}}{\t_{k+1}} (h(N_k)-h(\Nkp1)) & - & \sum_{k\in {\cal I}_2\cup {\cal I}_3}\frac{1-\t_{k+1}}{\t_{k+1}} (h(N_k)-h(\Nkp1))  \\
= \sum_{k=0}^{\bar k-1} \left( f_{N_k}(x_k)-f_{\Nkp1}(x_{k+1})\right ) & + &  \sum_{k\in  {\cal I}_{1}}  \frac{1-\t_{k+1}}{\t_{k+1}} (h(N_k)-h(\Nkp1))  \\
\le \sum_{k=0}^{\bar k-1} \left( f_{N_k}(x_k)-f_{\Nkp1}(x_{k+1})\right ) & + &  
\frac{1-\underline \t}{\underline \t} \sum_{k=0}^{ \bar k-1}h(N_k).
\end{eqnarray*}

This implies 
\begin{eqnarray}
&& \sum_{k=0}^{\bar k-1}\left( f_{N_k}(x_k)-f_{\Nkp1}(x_{k+1})\right )   +
\frac{1-\underline \t}{\underline \t} \sum_{k=0}^{ \bar k-1}h(N_k)  \ge \nonumber \\
&& \sum_{k\in {\cal I}_2\cup {\cal I}_3}\left( f_{N}(x_k)-f_{\Nkp1}(x_{k+1}) \right)  \label{bound_sum_f} \label{bound_sum_f}
\end{eqnarray}

Then, \eqref{bound_sum_f},  \eqref{sum_h}, \eqref{ared1} and 
$h(N_{\bar k})=0$ yield

\begin{eqnarray*}
 & & f_{N_0}(x_0)-f_{low}+ \left(\sigma\eta+\frac{1-\underline \t}{\underline \t}\right)  \, \frac{\xi-\underline{\t}f_{low}}{\eta^2(1-r)}   \\
  & &\quad  \ge \sum_{k=0}^{\bar k-1} \left( f_{N_k}(x_k)-f_{\Nkp1}(x_{k+1})  \right )+  
 \left(\sigma \eta+\frac{1-\underline \t}{\underline \t}\right) \sum_{k=0}^{\bar k-1}   h(N_k)\\
& &\quad \ge  \sum_{k\in {\cal I}_2\cup {\cal I}_3}\left( f_{N}(x_k)-f_{\Nkp1}(x_{k+1}) +\sigma \eta  h(\Nkp1) \right)\\
& &\quad \ge  (| {\cal I}_2|+ |{\cal I}_3|)  \frac{\eta}{2}    
\min \left\{\zeta_1\Gamma \varepsilon_g, \zeta_1 \sqrt{\frac{\z \underline{h} \z }{\mu}}, \Delta_0  \right\} \varepsilon_g,
 \end{eqnarray*}
 and claim  $ii)$ follows.   \epr

\vskip 5 pt

Considering that $\varepsilon_g$ is an optimality measure and $ \underline h$ is expected to be small, 
it is reasonable to suppose that 
\begin{equation}\label{delta0}
 \Delta_0 \ge \zeta_1 \max \left\{\Gamma \varepsilon_g,   \sqrt{\frac{\z \underline{h} \z }{\mu}}\right\}.
\end{equation}
Under this condition, Theorem \ref{th_complexity} gives the iteration complexity 
$$
|{\cal{I}}_1|+|{\cal{I}}_2|+ |{\cal{I}}_3|={\cal{O}} \left(\underline h^{-1}+\max\{ \varepsilon_g^{-2}, \, \underline h^{-\frac 1 2 }\varepsilon_g^{-1}\} \right).
$$
As a consequence,  for suitable values of $\underline h$, the worst-case iteration complexity $ {\cal{O}}(\varepsilon_g^{-2}) $
of the standard trust-region method  is retained, despite inaccuracy  in functions and gradients. This result is stated below,  where we count the number of iterations needed
to satisfy  $\|\nabla f_N(x_k)\|\le \varepsilon_g$ or $\|\nabla f_{N_{k+1}}(x_k)\|\le \varepsilon_g$ and $N_k=N$, i.e.,  iterations in ${\cal I}_1\cup{\cal I}_2\cup {\cal I}_3$ and iteration $\bar k$.
\begin{corollary} \label{corollary_compl}
Let Assumptions   \ref{assh}--\ref{Bk} hold. Assume furthermore $\gamma=1$ in  Algorithm {\sc iretr}.
Then,  there exists a constant $\nu_4>0$ such that Algorithm {\sc iretr} needs at most 
$$
\lceil \nu_4  \varepsilon_g^{-2}\rceil +1
$$
iterations, provided that  $\underline h^{-1}={\cal{O}}(\varepsilon_g^{-2} )$ { and \eqref{delta0} holds}.  
\end{corollary}
\vskip 5pt

In case $h(M)=(N-M)/N$, it holds $\underline h=1/N$ and $\underline h^{-1}={\cal{O}}(\varepsilon_g^{-2} )$    implies   $N= {\cal{O}}( \varepsilon_g^{-2})$.
In case  $N$ is larger, the number of iterations taken before full-accuracy is reached may deteriorate the complexity of the
standard trust-region approach.  

In order to derive the worst-case function evaluation complexity 
we need to bound the total number of trust-region reductions as each trust-region reduction 
calls for  one (possibly subsampled) function  evaluation at  trial point $x_k+p_k$.

\begin{theorem} \label{complexity_f}
Let Assumptions   \ref{assh}--\ref{Bk} hold.  Assume furthermore $\gamma=1$ in  Algorithm {\sc iretr}   and 
let ${\cal T}_j$ be the number of trust-region reductions at a generic iteration $j$ of the algorithm.
Then, for any $k\ge 1$,
$$
{
\sum_{j=0}^{k} {\cal T}_j \le \left\lceil \frac{\log(\underline \Delta /\Delta_0)}{\log(\zeta_1)}-k\frac{\log(\zeta_2)}{\log(\zeta_1)} \right\rceil }, 
$$
where 
$$
\underline{\Delta}=\min \left\{ \zeta_1  \sqrt{\frac{\eta (1-\eta)}{\kT +\mu} (1-r) \underline{h}},\zeta_1 \sqrt{\frac{\z \underline{ h } \z }{\mu}},\zeta_1 \Gamma \varepsilon_g ,\Delta_0 \right\}.
$$
\end{theorem}
\bpr
Let us proceed by induction. 
By the updating rules of the trust-region radius in Step 5 of Algorithm {\sc iretr}, {at termination of the  iteration $j=0$} we   have
$$
\Delta_1\in [ \zeta_1^{{\cal T}_0} \Delta_0, \zeta_2 \zeta_1^{{\cal T}_0} \Delta_0].
$$
Then, assume that at iteration $k\ge 1$ 
\begin{equation}\label{induction}
\Delta_k\in  [\zeta_1^{w_k} \Delta_0, \zeta_2^k \zeta_1^{w_k} \Delta_0],
\end{equation}
with $w_k=\sum_{j=0}^{k-1} {{\cal T}_j}$. At the end of iteration $k$, after ${\cal T}_k$ reductions of the trust-region radius  we have 
$$
\Delta_{k+1} \in [ \zeta_1^{{\cal T}_k} \Delta_k, \zeta_2 \zeta_1^{{\cal T}_k} \Delta_k ],
$$
and consequently, 
$$
\Delta_{k+1} \in [ \zeta_1^{w_{k+1}} \Delta_0, \zeta_2^{k+1} \zeta_1^{w_{k+1}} \Delta_0 ],
$$
i.e.,  \eqref{induction} holds for any $k\ge 1$. 
Taking into account that  Lemma \ref{lemma_Delta} ensures that iteration $k$ terminates with  $\Delta_k\ge \underline \Delta$,  
in  the adverse case 
where the initial $\Delta_k$ is given by $\zeta_2^k \zeta_1^{w_k} \Delta_0$ (see \eqref{induction}), at termination of iteration $k$ we are ensured that 
$$
\zeta_2^k \zeta_1^{w_{k+1}} \Delta_0\ge \underline \Delta.
$$
This yields the thesis, taking into account that $\zeta_1<1$.  $\Box$
\vskip 5 pt
Using the previous results we can now state our function evaluation complexity result.

\begin{corollary}
Let Assumptions   \ref{assh}--\ref{Bk} hold. Assume furthermore $\gamma=1$  in  Algorithm {\sc iretr}.  
Then,   if $\underline h^{-1}={\cal{O}}(\varepsilon_g^{-2}) $ and $\Delta_0$ satisfies \eqref{delta0}  and it is  independently of $\varepsilon_g$, there exists a constant $\nu_5$ such that Algorithm {\sc iretr} needs at most 
$$
\left\lceil    \nu_4\varepsilon_g^{-2} \left( 1- \frac{\log(\zeta_2)}{\log(\zeta_1)} \right)   -\frac{\log(\nu_5\varepsilon_g^{-1})}{\log(\zeta_1)}  \right\rceil
$$
function evaluations, where $\nu_4$ is  given in  Corollary \ref{corollary_compl}.

\end{corollary}
\bpr
Assumption $\underline h^{-1}={\cal{O}}(\varepsilon_g^{-2}) $, \eqref{delta0}  and $\Delta_0$ independent of $\varepsilon_g$ ensure
$\underline{\Delta}=\nu_5\varepsilon_g$, for some positive $\nu_5$. Then Corollary \ref{corollary_compl} and Theorem \ref{complexity_f} yield the thesis.  $\Box$

\section{Numerical experiments}\label{exps}
In this section we report on our numerical experience with Algorithm {\sc iretr} employing  the second order model (\ref{model}) 
and $D_{k+1}$  equal to a fixed fraction of $N_{k+1}$. 
Our aim is to show that our adaptive and deterministic strategy for 
choosing the  sample size $N_k$ and the use  of subsampled  functions, gradients and Hessians 
is effective and provides a gain in the overall computational cost with respect to a standard trust-region approach. 
To this end, we   compare our method with  ``standard'' trust-region implementations, i.e. implementations where functions and gradients are computed at full accuracy too.
 Specifically, we compare with the implementation, named {\sc statr\_sh},  employing full   functions and gradients and subsampled Hessian 
$B_k$  as  in (\ref{Bk_sampling}) with  $\Dkp1=\left\lceil 0.1N \right \rceil$, and with   the implementation, named   {\sc statr\_fh}, where functions, first and second order derivatives  are computed at full accuracy. 

All the results have been obtained running    a Matlab  R2019b code  on an  Intel Core i5-6600K CPU 3.50 GHz x 4,  16.0GB RAM.

\subsection{Test problems}
%
We tested our method both on convex and nonconvex problems arising in binary classification problems.
Let $\{(a_i, b_i)\}_{i=1}^N$ denote the pairs forming the data set with  $a_i \in \IR^n$ 
being  the vector containing the entries of the $i$-th example   and $b_i$ being its label. 
The data set we employed   are displayed in Table  \ref{test}. 
In the table for each data set we report the number $N$ of training examples  
and the dimension $n$ of each instance. Moreover we report   the number of elements in the testing set $N_T$.

We performed a logistic regression to solve classification problems associated to 
the data sets  {\sc Mushrooms},  {\sc Cina0} and {\sc Gisette}. In this case   $b_i\in \{-1, +1\}$ and 
the  strongly convex objective function is given by  the logistic loss with $\ell_2$-regularization
$$
f_N(x)=\frac 1 N\sum_{i=1}^N\log(1+e^{-b_i a_i^Tx})+\frac{1}{2N} \|x\|^2.
$$

Classification problems associated with the remaining  data  sets were solved using the sigmoid function and  least-squares loss. 
Here $b_i\in \{0, +1\}$ and  the non-convex objective function has the form
$$
f_N(x)=\frac 1 N\sum_{i=1}^N \left(b_i-\frac{1}{1+e^{-a_i^Tx}} \right)^2.
$$

\begin{small}
\begin{table} 
\begin{center}
\begin{tabular}{l  |r r |r }
&     \multicolumn{2}{r|} { Training set }& \multicolumn{1}{r} { Testing set }\\
\hline 
Data set & $N$ & $n$  &  $N_T$   \\  \hline
{\sc Mushrooms } \cite{UCI}& 5000 & 112 &  3124\\
{\sc Cina0}  \cite{CINA} &  10000 & 132 & 6033 \\
{\sc Gisette} \cite{UCI}& 5000 & 5000 & 1000 \\
{\sc A9a} \cite{UCI}& 22793& 123& 9768\\
 {\sc Covertype} \cite{UCI}& 464810& 54 & 116202 \\
{\sc Ijcnn1} \cite{libsvm}&49990 & 22& 91701\\
{\sc Mnist}  \cite{mnist}& 60000& 784& 10000  \\
{\sc Htru2}\cite{UCI} &10000  & 8  & 7898  \\
\hline
\end{tabular}
\caption{Data sets used}\label{test}
\end{center}
\end{table}
\end{small}

\subsection{Implementation issues}
The trust-region parameters of  the procedures under comparison are fixed as 
$$
\quad \Delta_0=10,   \quad \tau=0.1,  \quad \eta=0.1, \quad \zeta_1=0.5, \quad \zeta_2=1.2. 
$$
The trust-region problem is solved approximately using CG-Steihaug method \cite{cgt}.
The Conjugate Gradient (CG) method is applied without preconditioning and the procedure is stopped when the relative residual 
becomes smaller than $10^{-3}$ or a maximum of $100$ iterations is performed. 
In Step 5, in case of successful   iterations,  we update the trust-region radius as follows. If
 $\ared(\t_{k+1})/\pred(\t_{k+1})\ge 1.1$ we set $\Delta_{k+1}^{(0)}=\zeta_2 \Delta_k^{({\cal T}_k)}$, otherwise we  set  $\Delta_{k+1}^{(0)}= \Delta_k^{({\cal T}_k)}$.
 \medskip \noindent

Focusing on Algorithms {\sc iretr},   we tested two rules for choosing the sample size. 
In the first implementation, later referred to as {\sc iretr\_d}, the sample size varies  dynamically.  
The infeasibility measure $h$  and the initialization  parameters for inexact restoration are: 
$$
h(M)=\frac{N-M}{N}, \quad N_0=\left \lceil 0.1\,N \right \rceil, \quad  \t_0=0.9.
$$
The parameters $\gamma=1, \, \mu=100/N$ are used in  \eqref{new2}.
The updating rules for choosing $\Ntilde$, $\Nkp1$ in Steps 1 and 2  are the following: 
\begin{eqnarray*}
\Ntilde&=&\min\{N, \left\lceil 1.2\, \Nk\right\rceil\},\\
\Nkp1&=&\left\{\begin{array}{ll}       
\left\lceil  \Ntilde-10^2 \Delta_k^{1+\gamma} \right\rceil   & \mbox{ if } \left\lceil \Ntilde-10^2\Delta_k^{1+\gamma} \right\rceil \in [N_0, 0.95 N],\\
&   \\
\Ntilde & \mbox{ if } \left\lceil \Ntilde-10^2\Delta_k^{1+\gamma} \right\rceil <N_0, \\
&   \\
N & \mbox{ if } \left\lceil \Ntilde-10^2\Delta_k^{1+\gamma} \right\rceil >0.95 N.
\end{array}\right.
\end{eqnarray*}
We note that the choice of $\Ntilde$ falls into (\ref{feas}) with $r=(N-0.2)/N$. 

In the second implementation, we set again 
$$
h(M)=\frac{N-M}{N}, \quad  \t_0=0.9.
$$
Then, the sample size $\Nkp1$ is increased according the geometric growth:  
$$
N_0=\left \lceil 0.1\,N\right \rceil,\qquad \Nkp1=\Ntilde = \min\{N, \left\lceil 1.2\, \Nk\right\rceil\}.
$$
We will refer to this implementation as  {\sc iretr\_gg}. We note that this choice of $\Nkp1$ amount to choosing $\mu=0$ in (\ref{new2}).

In both implementations {\sc iretr\_d} and {\sc iretr\_gg}
the first time that $\Nk=\Nkp1=N$ occurs, then the value of the trust-region radius is set to  $\Delta_k^{({\cal T}_k)}=\max\{1,\,\Delta_k^{({\cal T}_k)} \}$.
Moreover, the Hessian   matrix $B_k$ is formed via  \eqref{Bk_sampling} with
$$
\Dkp1=\left\lceil 0.1\,  \Nkp1\right\rceil, \ \  \forall k\ge 0.
$$
Thus, the Hessian sample size changes dynamically until the full sample for function and gradient is  reached. 
The sets $I_{\Nkp1}$ and $I_{\Dkp1}$ are generated using the {\tt Matlab} function {\tt randsample} with no replacement.
When  the sample size $\Nkp1$ is increased, the new sample set can be computed from scratch or can be obtained randomly adding new samples to the previous sample set. Despite this latter choice produces  computational savings, in view of a truly random process  we generate each $I_{\Nkp1}$ from scratch.

Concerning  the stopping criteria, for all the algorithms under comparison, we imposed a maximum of  $1000 $ iterations
and we declared a successful termination  when one of the two following conditions is met 
\begin{equation}\label{stop}
\|\nabla f_{\Nk}(x_k)\|\le  \varphi  , \qquad |f_{\Nk}(x_k)- f_{N_{k-1}}(x_{k-1})  | \le \varphi |f_{\Nk}(x_k)|,
\end{equation}
with $\varphi=10^{-4}$.
We underline that for {\sc iretr\_d} and {\sc iretr\_gg} the above checks are on possibly subsampled functions and gradients and we allow 
for termination before full precision is reached. 

The initial guess is $x_0=(0,\ldots, 0)^T$ for all runs.

\subsection{Numerical  results}
The first set of  results presented shows  the performance of  Algorithms {\sc iretr\_d}, {\sc iretr\_gg}, {\sc statr\_sh} and  {\sc statr\_fh}.  
In our test problems, the main cost in the computation of  $\phi_i$ for any $1\le i\le N$  is the scalar product $a_i^Tx$. Once this product is evaluated,
it can be reused for computing $\nabla \phi_i$ and  $\nabla^2 \phi_i$.  In particular, 
computing $\nabla^2 \phi_i$ times a vector  $v$  at each CG iteration requires a scalar product $a_i^Tv$ i.e., it is as expensive as
evaluating $\phi_i$. Therefore, if one full function evaluation is denoted as   {\tt nfe}, 
computing $f_M$  costs $\displaystyle \frac{M}{N}${\tt nfe} while each CG iteration costs $\displaystyle \frac{\Dkp1}{N}${\tt nfe}. 
Since the selection of sets  $I_{\Nkp1}$ and $I_{\Dkp1}$ in Algorithms  {\sc iretr\_d}, {\sc iretr\_gg} and {\sc statr\_sh}  is random, the cost associated to such  algorithms is  measured on average over 50 runs.

In Table \ref{table_result}  for each method and for each data set we report the number  {\tt nfe} of full function evaluations 
performed  and the percentage of saving obtained  by   Algorithm {\sc iretr\_d}  with respect to   {\sc iretr\_gg},
{\sc statr\_sh}  and to {\sc statr\_fh}. 
First, we can observe that Algorithm {\sc iretr\_d}  is in general less costly than  the variant {\sc iretr\_gg}; this  indicates
that the dynamic choice of the sample size, aiming to make slow progress to full precision, is effective and does not deteriorate the performance of 
{\sc iretr} when the geometrical growth of the sample size is the most effective (see the results for {\sc Htru2}).
Second, we observe a remarkable  saving of both {\sc iretr\_d}  and   {\sc iretr\_gg} with respect to the full standard trust-region for all the data sets used; compared to {\sc statr\_sh}
the saving   is lower, as expected,  but still considerable overall.

\begin{small}
\begin{table}   
\begin{center}
\begin{tabular}{l  |c |rrr}
Data set & {\tt nfe} & \multicolumn{3}{|c } { {\tt nfe(save)} } \\  \hline
& {\sc iretr\_d}  &  {\sc iretr\_gg} &{\sc statr\_sh} &  {\sc statr\_fh} \\  \hline
{\sc Mushrooms} & 27 & 30 (10\%) & 51 (47\%) &108 (75\%)    \\
{\sc Cina0}  &  88 & 99 (11\%) & 96 (\z 8\%)& 416 (78\%)\\
{\sc Gisette} &346& 362 (\z4\%)& 432 (20\%) & 594 (42\%) \\
{\sc A9a} &  22 & 25 (12\%) &45 (51\%)  & 445 (95\%)  \\
 {\sc Covertype} & 17 & 23 (26\%) &  48 (65\%)& 698 (98\%)  \\
{\sc Ijcnn1} &20 & 25 (20\%) &36 (44\%) & 128 (84\%)  \\
{\sc  Mnist  } & 46 & 50 (\z8\%) & 58 (20\%)  & 955 (95\%)  \\
{\sc Htru2} &38 & 37 (\ -3\%)  & 43 (12\%) & 87 (56\%)\\  
\hline
\end{tabular}
\caption{Function evaluations performed by    {\sc iretr\_d}, {\sc iretr\_gg}, {\sc statr\_sh}  and   {\sc statr\_fh}  and saving obtained  by  
{\sc iretr\_d} over {\sc iretr\_gg}, {\sc statr\_sh}  and   {\sc statr\_fh}. }   \label{table_result}
\end{center}
\end{table}
\end{small} 
%
%
%
%
%

\begin{figure}
\centering
  \includegraphics[height=0.3\textheight,width=1\textwidth]{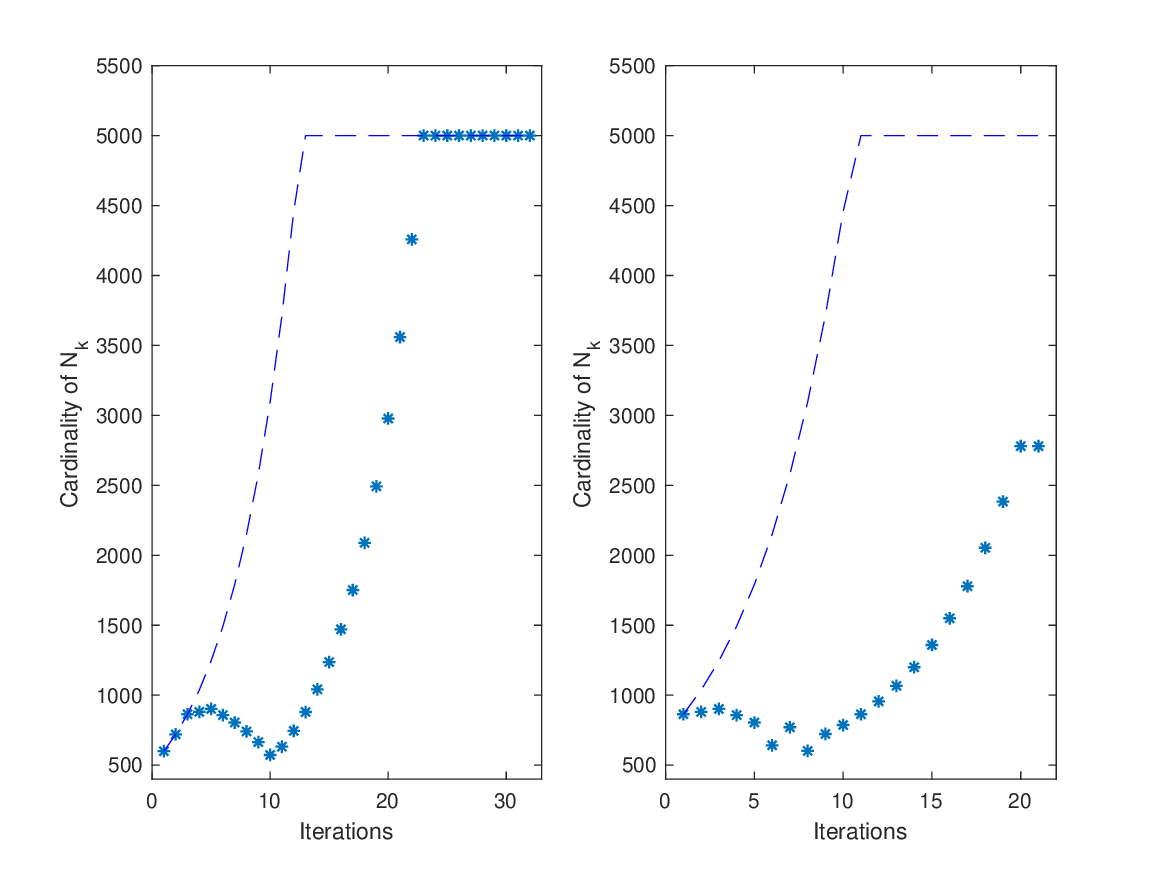}
  \hspace*{-0.5cm}
  \caption{{\sc Mushrooms} data set. $\Nk$ versus {\sc iretr\_d} iterations (`` * ''), sample size  $\Nkp1=(1.2)^k N_0$ (dashed line).} 
  \label{fignkM}
\end{figure}

\begin{figure}
\centering
  \includegraphics[height=0.3\textheight,width=1\textwidth]{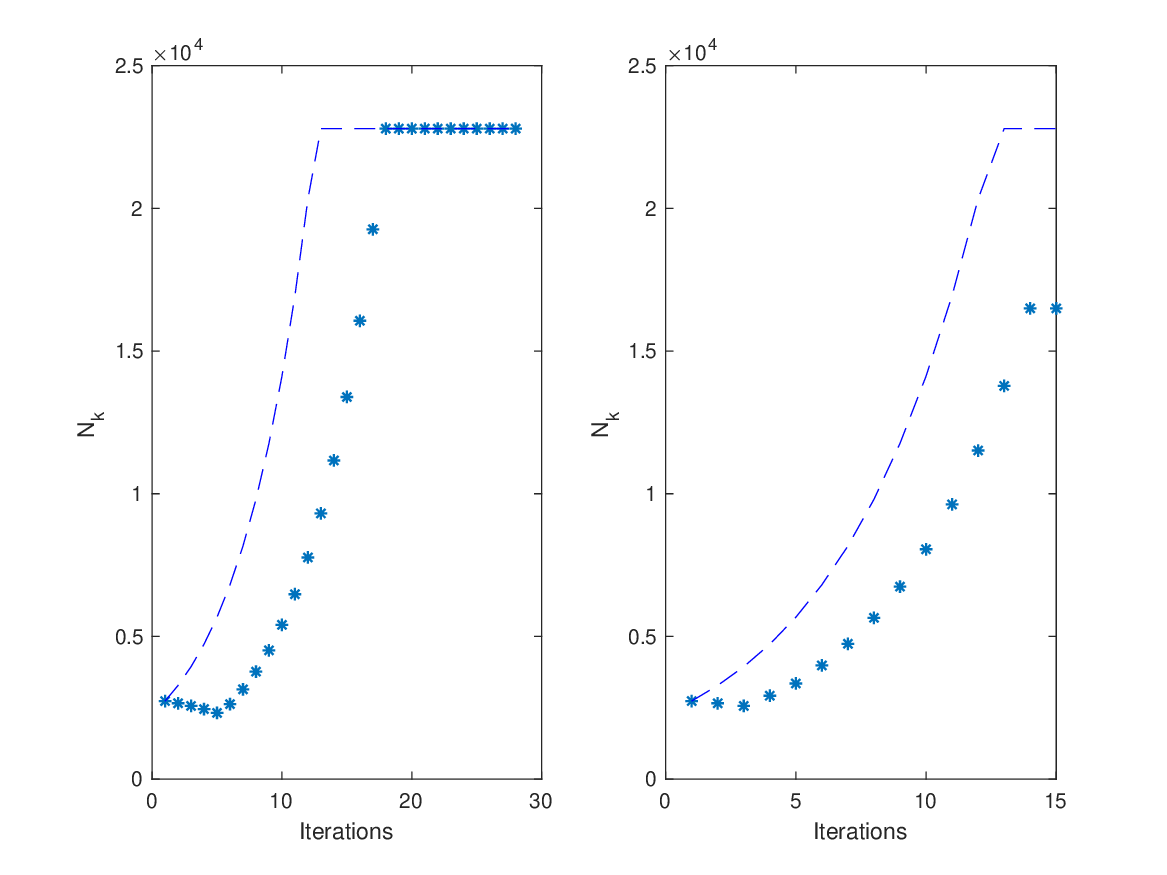}
  \caption{{\sc A9a} data set. $\Nk$ versus {\sc iretr\_d} iterations (`` * ''), sample size  $\Nkp1=(1.2)^k N_0$ (dashed line).} 
  \label{fignkA}
\end{figure}

\begin{figure}
\centering
  \includegraphics[height=0.3\textheight,width=1\textwidth]{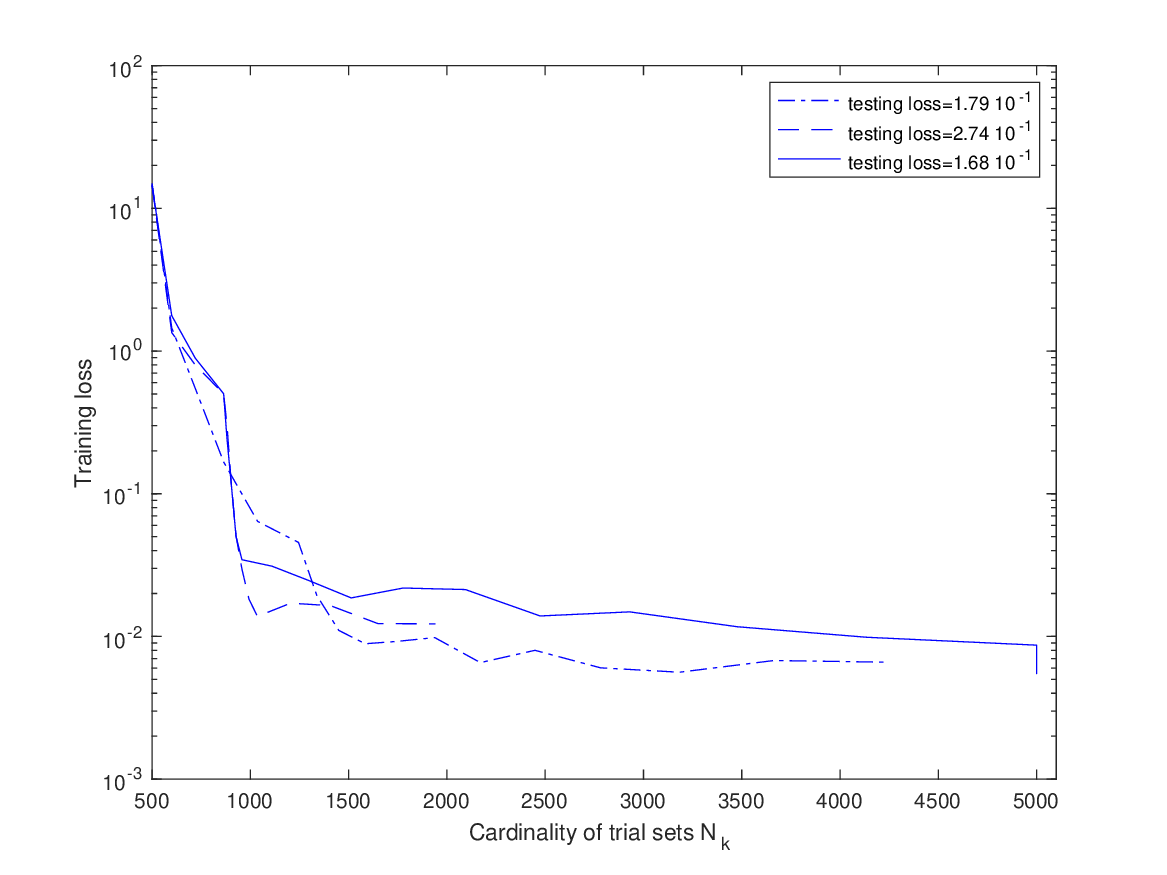}
  \caption{  {\sc Mushrooms} data set, $N$=5000. Training loss  versus   $\Nk$ and testing loss at termination using {\sc iretr\_d}. Values of $\Nk$ at termination: 
	1941 (dashed  line); 4241 (dash-dotted line);  5000 (solid line).}
  \label{fibjcnn1}
\end{figure}

To give more insight into the two implementations {\sc iretr\_d},
in Figures  \ref{fignkM} and  \ref{fignkA} we  plot the sample size $N_k$ versus the iterations for {\sc Mushrooms} and   {\sc A9a} problems.
The dashed line plots $N_{k+1}=\left \lceil (1.2)^kN_0\right \rceil$ versus iterations, that is the sample size corresponding to the geometric growth used in {\sc iretr\_gg}.
The increase of $N_k$ along iterations  in {\sc iretr\_d}   is  considerably slower 
than that provided by the geometric growth; 
in two runs,  the cardinality  $\Nk$ in   {\sc iretr\_d}  reaches the value $N$, as expected from the theory, but in the
first phase of the iterative process it is a small fraction of $N$ and decreases at some iterations. 
In the other two runs, {\sc iretr\_d}  does not reach full precision, iterations terminate  
with a cardinality $\Nkp1=2780$, corresponding to the 56\% of the training set and $\Nkp1=16495$,  corresponding to the 72\% of the training set, respectively.
In fact, despite the adaptive strategy of  {\sc iretr}  yields   $N_k=N$ for $k$ sufficiently large, 
our stopping rule (\ref{stop}) is applied on possibly subsampled functions and gradients.   
This feature is in accordance with the motivations for using subsampling: data in a training set show redundancy and in general using
subsets of the sample data is enough to provide a small testing error.
At this regard, consider Figure \ref{fibjcnn1}
related to  the data set  {\sc Mushrooms}, $N=5000$.  At each iteration and for   three runs corresponding 
to different sample sizes at termination,  we plot the training loss $f_{\Nk}(x_k)$ versus the value of $N_k$;  at termination:
$N_k$ =1941 (dashed  line),  $N_k$= 4241 (dash-dotted line), $N_k=N$ (solid line).  We also display the
testing loss $f_{N_T}$ at termination.   Although in two runs the final sample size is approximately 
39\% and 85\%  of the data in the training set,  interestingly 
the testing loss is in between   $1\cdot 10^{-1}$ and  $3\cdot 10^{-1}$  in all runs. Thus,  monitoring the values 
of subsampled  functions and gradients  in (\ref{stop}) is effective.
 
\begin{figure}

\centering
  \includegraphics[height=0.4\textheight,width=1\textwidth]{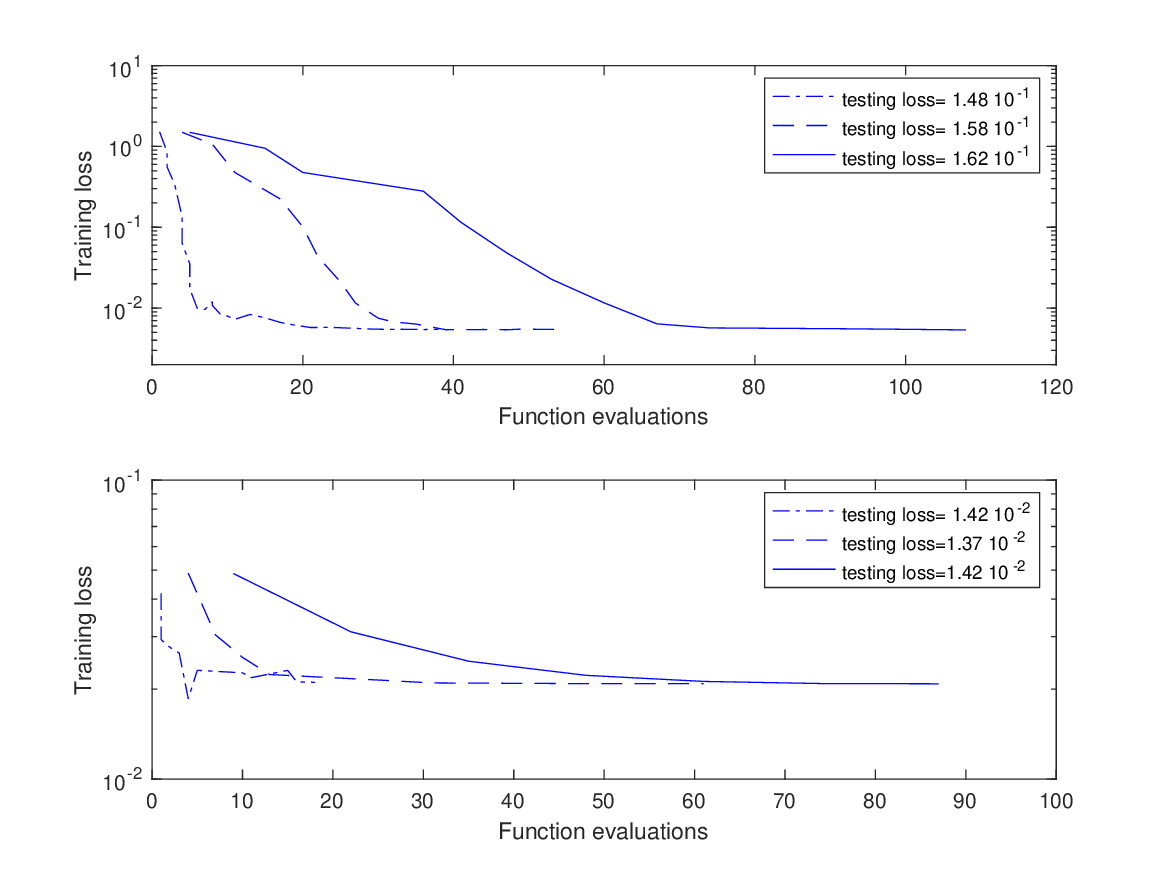}
  \caption{Training loss versus function evaluations and testing loss:  {\sc iretr\_d} (dash-dotted line);
{\sc statr\_sh}    (dashed line) and   {\sc statr\_fh} (solid line). Upper: {\sc Mushrooms} data set, lower: {\sc Htru2} data set.} 
  \label{figmushpulsar}
\end{figure}

The previous discussion is supported by further observations.
In Figure  \ref{figmushpulsar}, we plot the  value  of the training loss versus the number of function evaluations required to 
solve {\sc Mushrooms}  and {\sc Htru2} problems with   {\sc iretr\_d}, {\sc statr\_sh}    and   {\sc statr\_fh}.
In these runs, {\sc iretr\_d} terminates with $\Nk=N$ in {\sc Mushrooms}   problem while terminates with $\Nk=7426$ (74\% of the samples)  in {\sc Htru2} problem.
At termination, the values of both the training loss and the testing loss  provided by the three methods   are similar and this feature further supports  both
termination before full precision is reached and the
inexact restoration approach  for handling  subsampled functions and derivatives.

\begin{figure}
\centering
  \includegraphics[height=0.3\textheight,width=1\textwidth]{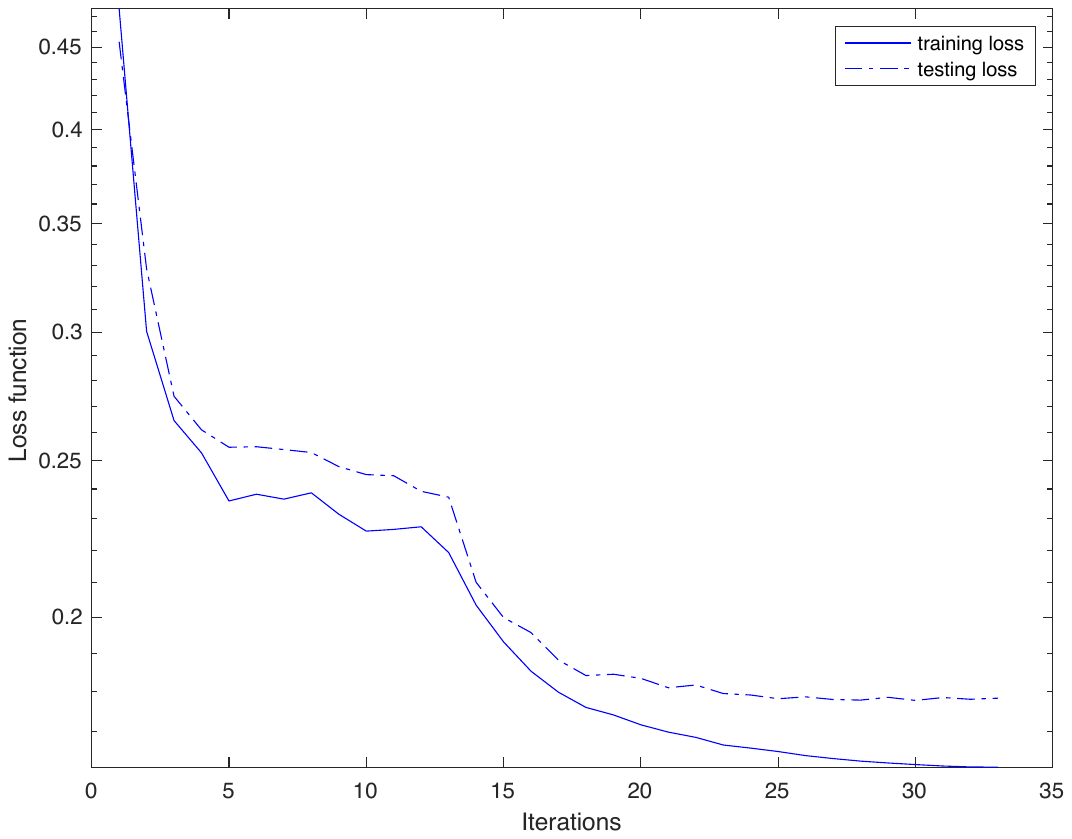}
  \caption{ {\sc Cina0} data set: training and testing  loss function versus iterations computed by {\sc itetr\_d}. Stopping threshold $\varphi=10^{-8}$.} 
  \label{figcina}
\end{figure}

Finally, Figure   \ref{figcina} refers to  the dataset {\sc Cina0} and displays the values of the training and testing logistic loss 
along the iterations  of {\sc iretr\_d} using the tolerance  $\varphi=10^{-8}$ in  \eqref{stop}. In the progress of the iterations  the loss
values  settle  and performing   the last thirteen   iterations  is pointless.

%
%

\noindent{\bf Acknowledgement} Dedicated with friendship to Jos\'e Mario Mart\'{\i}nez for his outstanding scientific contributions.

\end{document}